\UseRawInputEncoding
\documentclass[12pt,reqno]{amsart}
\usepackage{caption,amsmath,amssymb,amsfonts,amsthm,latexsym,graphicx,multirow,color,hyperref,enumerate}

\newtheorem{lem}{Lemma}[section]%
\newtheorem{theorem}[lem]{Theorem}%
\newtheorem{exam}[lem]{Example}%
\newtheorem{prob}[lem]{Problem}%
\newtheorem{prop}[lem]{Proposition}%

\def\a{\alpha} \def\b{\beta} \def\g{\gamma} \def\d{\delta} 
 \def\s{\sigma}   
  
 \def\L{\Lambda}
 
\def\G{\Gamma}

 \def\lg{\langle} \def\rg{\rangle}
\def\nd{\mathrel{\bigm|\kern-.7em/}}

\def\f{\noindent}

\def\PSL{\hbox{\rm PSL}}
\def\AGL{\hbox{\rm AGL}}
\def\PL{\hbox{\rm P$\Gamma$L}}
\def\AGammaL{\hbox{\rm A$\Gamma$L}}
\def\SL{\hbox{\rm SL}}
\def\GL{\hbox{\rm GL}}
\def\GF{\hbox{\rm GF}}
\def\PGL{\hbox{\rm PGL}}
\def\PSU{\hbox{\rm PSU}}

\def\Aut{\hbox{\rm Aut}}

\def\soc{\hbox{\rm soc}}

\def\BiCay{\hbox{\rm BiCay}}
\def\Cay{\hbox{\rm Cay}}

\def\mod{\hbox{\rm mod }}

\def\S{\hbox{\rm S}}
\def\F{\hbox{\rm Frob}}
\def\demo{\f {\bf Proof.}\hskip10pt}

\renewcommand\qed{\hskip10pt $\Box$\vspace{3mm}}
\def\mz{{\mathbb Z}}
\def\C{{C}}

\def\A{\hbox{\rm A}}
\def\H{\mathcal{H}}
\def\M{\mathcal{M}}
\def\N{\mathcal{N}}
\def\R{\mathcal{R}}
\def\T{\mathcal{T}}
\def\L{\mathcal{L}}

\def\K{{\bf K}}
\def\1{{\bf 1}}

\def\Ome{{\rm\Omega}}

\begin{document}

\title[$3$-CSH locally $2\K_n$ graphs and $s$-arc-transitive graphs ]
{Finite $3$-connected-set-homogeneous locally $2\K_n$ graphs and $s$-arc-transitive graphs }
\author{Jin-Xin Zhou}%
\address{School of Mathematics and Statistics, Beijing Jiaotong University, Beijing, 100044, China}
 \email{jxzhou@bjtu.edu.cn} %

\thanks{2010 Mathematics Classification: 05C25, 05E18, 20B25.}



\maketitle

\begin{abstract}
In this paper, all graphs are assumed to be finite. For $s\geq 1$ and a graph $\G$, if for every pair of isomorphic connected induced subgraphs on at most $s$ vertices there exists an automorphism of $\G$ mapping the first to the second, then we say that $\G$ is $s$-connected-set-homogeneous, and if every isomorphism between two isomorphic connected induced subgraphs on at most $s$ vertices can be extended to an automorphism of $\G$, then we say that $\G$ is $s$-connected-homogeneous. For $n\geq 1$, a graph $\G$ is said to be locally $2\K_n$ if the subgraph $[\G(u)]$ induced on the set of vertices of $\G$ adjacent to a given vertex $u$ is isomorphic to $2\K_n$.

Note that $2$-connected-set-homogeneous but not $2$-connected-homogeneous graphs are just the half-arc-transitive graphs which are a quite active topic in algebraic graph theory. Motivated by this, we posed the problem of characterizing or classifying $3$-connected-set-homogeneous graphs of girth $3$ which are not $3$-connected-homogeneous in (Eur. J. Combin. 93 (2021) 103275). Until now, there have been only two known families of $3$-connected-set-homogeneous graphs of girth $3$ which are not $3$-connected-homogeneous, and these graphs are locally $2\K_n$ with $n=2$ or $4$. In this paper, we complete the classification of finite $3$-connected-set-homogeneous graphs which are locally $2\K_n$ with $n\geq 2$, and all such graphs are line graphs of some specific $2$-arc-transitive graphs. Furthermore, we give a good description of finite $3$-connected-set-homogeneous but not $3$-connected-homogeneous graphs which are locally $2\K_n$ and have solvable automorphism groups. This is then used to construct some new $3$-connected-set-homogeneous but not $3$-connected-homogeneous graphs as well as some new $2$-arc-transitive graphs. 

\bigskip

\end{abstract}
\thispagestyle{empty}

\section{Introduction}
The main purpose of this paper is to give a partial answer to a problem raised by the author in \cite{Zhou-EJC-2021} about the $3$-connected-set-homogeneous graphs. As a by-product of this investigation, we disprove a conjecture posed by Feng and Kwak in their 2006 paper on trivalent symmetric graphs of order twice a prime power~\cite{Feng-Kwak}, and answer a question on $3$-arc-transitive graphs posed  by Li, Seress and Song in \cite{Li-S-Song-2015-JA}, and we also correct an error in \cite{LiLingMa} about tetravalent $3$-arc-regular Cayley graphs. Before proceeding, we give some background to this topic, and set some notation.

A graph is said to be {\em regular} if each of its vertices is adjacent to $k$ vertices for some constant positive integer $k$.
Let $\G$ be a graph. We use $V(\G)$, $E(\G)$ and $\Aut(\G)$ to denote its vertex set, edge set and full automorphism group, respectively. For $B\subseteq V(\G)$, $[B]$ denotes the subgraph induced by $B$. For a vertex $v$ of $\G$, let $\G(v)$ be the set of vertices adjacent to $v$. Let $G\leq \Aut(\G)$. Denote by $G_v$ the subgroup of $G$ fixing $v$, by $G_v^{\G(v)}$ the permutation subgroup on $\G(v)$ induced by $G_v$, and by $G_v^{[1]}$ the subgroup of $G_v$ fixing every vertex in $\G(v)$. For a positive integer $n$, we say that $\G$ is {\em locally $2\K_n$} if $[\G(v)]\cong 2\K_n$ for each $v\in V(\G)$, where $2\K_n$ means the disjoint union of two copies of $\K_n$.

For a positive integer $n$, denote by $C_n$ the cyclic group of order $n$, 
and by $\A_n$ and $\S_n$ the alternating group and symmetric group of degree $n$, respectively. For two groups $M$ and $N$, $N\rtimes M$ denotes a semidirect product of $N$ by $M$, and $N\wr M$ the wreath product of $N$ by $M$.
See Section~\ref{sec:2} for other unexplained terms.

For a positive integer $s$ and a graph $\G$, if for any pair of isomorphic connected induced subgraphs of $\G$ on at most $s$ vertices there is an automorphism of $\G$ mapping the first to the second, then we say that $\G$ is {\em $s$-connected-set-homogeneous}, or {\em $s$-CSH}, and if every isomorphism between two isomorphic connected induced subgraphs on at most $s$ vertices can be extended to an automorphism of $\G$, then we say that $\G$ is {\em $s$-connected-homogeneous}, or {\em $s$-CH}. A graph is said to be {\em connected-set-homogeneous} or {\em connected-homogeneous} if it is $s$-CSH or $s$-CH, respectively, for all positive integers $s$.

$s$-CSH or $s$-CH graphs have received a lot of attention in the literature. For example, in 1978, Gardiner~\cite{Gardiner1978} gave a classification of finite connected-set-homogeneous graphs, in 2009, Gray~\cite{Gray} classified infinite $3$-CSH or $3$-CH graphs with more than one end, and Devillers et al. \cite{DFPZ,Li-Zhou} investigated the finite $k$-CH graphs with $k\geq 3$. For more results related to $s$-CSH or $s$-CH graphs, we refer the reader to \cite{Droste,Enomoto,Gray-Macpherson,Lachlan-Woodrow}.

Clearly, a graph is $1$-CSH or $1$-CH if and only if it is vertex-transitive. Furthermore, $2$-CH graphs are precisely regular arc-transitive graphs, and every $2$-CSH graph is vertex- and edge-transitive. A graph is said to be {\em half-arc-transitive} if it is $2$-CSH but not $2$-CH. In 1966, Tutte~\cite{Tutte1966} initiated the study of half-arc-transitive graphs, and he proved that the valency of a half-arc-transitive graph must be even, and a few years latter, Bouwer~\cite{Bouwer} constructed the first family of half-arc-transitive graphs. Following this pioneering work, half-arc-transitive graphs have been extensively studied over the last half a century, and numerous papers have been published on this class of graphs (see, for example, the survey papers \cite{CPS,Marusic1998} and recent
papers \cite{Spiga,Spiga-Xia2021JCTA,XiaJCTB2021,Zhou-JACO}).

In this paper, we are interesting in $3$-CSH but not $3$-CH graphs which are a natural generalization of half-arc-transitive graphs.
Note that a graph of girth at least $4$ is $3$-CSH but not $3$-CH if and only if it is $2$-path-transitive but not $2$-arc-transitive.
In 1996, Conder and Praeger \cite{Conder-Praeger} initiated the study of 2-path-transitive graphs, and more than ten years later, Li and Zhang~\cite{LZ2,LZ} systematically investigated $2$-path-transitive graphs which are not $2$-arc-transitive. Motivated by this, the author~\cite{Zhou-EJC-2021} began the study of $3$-CSH but not $3$-CH graphs of girth $3$, and we proved the existence of such graphs and proposed the following problem.

\begin{prob}{\rm \cite[Problem~B]{Zhou-EJC-2021}}\label{prob}
Characterize or classify $3$-connected-set-homogeneous graphs of girth $3$ which are not $3$-connected-homogeneous.
\end{prob}

In this paper, we shall partially solve this problem by classifying $3$-CSH graphs which are locally $2\K_n$ with $n\geq 2$. This was partially motivated by our previous work in \cite{Zhou-EJC-2021}, where we proved the existence of $3$-CSH but not $3$-CH graphs of girth $3$ by constructing some $3$-CSH graphs which are locally $2\K_n$ with $n=2$ or $4$. Our main results show that there is a close relationship between $3$-CSH graphs which are locally $2\K_n$ and $s$-arc-transitive graphs with $s\geq 2$. Constructing or classifying $s$-arc-transitive graphs with $s\geq 2$ has been a perennially active topic in the area of algebraic graph theory; see, for example, \cite{Li2001,Li-2005-PAMS,Li-S-Song-2015-JA,Praeger1993}. This is another motivation for us to study $3$-CSH graphs which are locally $2\K_n$. Before stating our main results, we introduce some terminology.

For $s\geq 0$, an {\em $s$-arc} in $\G$ is an ordered $(s+1)$-tuple $(v_0,v_1, \cdots ,v_{s-1},v_s)$ of vertices of $\G$ such that $v_{i-1}$ is adjacent to $v_i$ for $1\leq i\leq s$ and $v_{i-1} \neq v_{i+1}$ for $1 \leq i \leq s-1$. For some group $G$ of automorphisms of $\G$, we say that the graph $\G$ is {\em $(G, s)$-arc-transitive} if $\G$ is regular and $G$ is transitive on the set of $s$-arcs in $\G$; $\G$ is $(G,s)$-{\em arc-regular} if $G$ is regular on the set of $s$-arcs of $\G$. When $G=\Aut(\G)$, a $(G, s)$-arc-transitive or $(G, s)$-arc-regular graph $\G$ is simply called {\em $s$-arc-transitive} or {\em $s$-arc-regular}, respectively. 

The line graph $L(\G)$ of a graph $\G$ is the graph whose vertices are the edges of $\G$, with two edges adjacent in $L(\G)$ if they have a vertex in common. A graph $\G$ is said to be {\em locally $3$-transitive} if the vertex-stabilizer $\Aut(\G)_v$ of $v\in V(\G)$ acts $3$-transitively on $\G(v)$.

Now we state our first main theorem.

\begin{theorem}\label{mainth}
Let $\G$ be a locally $2\K_n$ graph with $n\geq 2$. Then $\G$ is $3$-connected-homogeneous if and only if $\G$ is isomorphic to the line graph $L(\Sigma)$ of a $3$-arc-transitive and locally $3$-transitive graph $\Sigma$.
\end{theorem}

By Theorem~\ref{mainth}, to construct $3$-CH locally $2\K_n$ graphs with $n\geq 2$, it is equivalent to construct $3$-arc-transitive and locally $3$-transitive graphs. The following theorem characterizes vertex stabilizers of $3$-arc-transitive graphs.

\begin{theorem}{\rm \cite[Theorem~4.2]{Li-S-Song-2015-JA}}\label{th:3-arc-tran-stabi}
For a $(G,3)$-arc-transitive graph $\G$ of valency $k$ and a $2$-arc $(w,u,v)$,  at least one of the following holds:
  \begin{enumerate}
  \item [{\rm (i)}]\ $G_u^{[1]}$ is transitive on $\G(w)-\{u\}$, or
  \item [{\rm (ii)}]\ $G_u=\A_7$ or $\S_7$, and $k=7$, or
  \item [{\rm (iii)}]\ $C_p^f\unlhd G_u^{\G(u)}\leq {\rm A\Gamma L}_1(p^f)$, the number of $G_u^{[1]}$-orbits on $\G(w)-\{u\}$ divides ${\rm gcd}(p^f-1, f)^2$, and $(G_u^{\G(u)})_{wv}\leq C_f$.
\end{enumerate}
\end{theorem}

There are infinitely many $3$-arc-transitive and locally $3$-transitive graphs satisfying the condition in Theorem~\ref{th:3-arc-tran-stabi}~(i). For example, the complete bipartite graph $\K_{n,n}$ is a $3$-arc-transitive and locally $3$-transitive graph for each $n\geq 3$. For more examples of $3$-arc-transitive and locally $3$-transitive graphs satisfying the condition in Theorem~\ref{th:3-arc-tran-stabi}~(i), we refer the reader to \cite{Li2001,LZ3}. It is easy to see that every $3$-arc-transitive graph satisfying the condition in Theorem~\ref{th:3-arc-tran-stabi}~(ii) is locally $3$-transitive. 
Recently, Giudici and King~\cite{Giudici-King-3-arc-transitive} gave a classification of edge-primitive $3$-arc-transitive graphs satisfying the condition in Theorem~\ref{th:3-arc-tran-stabi}~(ii), and there are two sporadic and eight infinite families of such graphs. 

In part (c) of \cite[Remarks on Theorem~4.2]{Li-S-Song-2015-JA}, Li et al. wrote  ``It is not known whether there are 3-arc-transitive graphs satisfying the condition in part (iii) of Theorem 4.2". It is easy to see that there does not exist a $3$-arc-transitive and locally $3$-transitive graph satisfying the condition in Theorem~\ref{th:3-arc-tran-stabi}~(iii). However, our next result shows that there do exist $(G,3)$-arc-transitive graphs satisfying the condition in Theorem~\ref{th:3-arc-tran-stabi}~(iii).

\begin{prop}\label{lem:sol-s-tran-stabi}
Let $p$ be a prime and $f$ be a positive integer. If $p^f-1$ is not coprime to $f$, then there exists $G\leq\Aut(\K_{p^f,p^f})$ such that $\K_{p^f,p^f}$ is a $(G,3)$-arc-transitive graphs satisfying the condition in Theorem~\ref{th:3-arc-tran-stabi}~(iii).
\end{prop}

\medskip
Following \cite{Morgan}, we say that a pentavalent symmetric graph $\G$ is of {\em type ${\mathcal{Q}_2^6}$} if $\Aut(\G)_u\cong \F(20)\times C_2$ and $\Aut(\G)_{\{u,v\}}\cong {\rm M}_{16}$, where $\{u,v\}$ is an edge of $\G$,
\[\begin{array}{l}
\F(20)=\lg a,b\mid a^{5}=b^4=1, b^{-1}ab=a^2\rg\ {\rm and}\ {\rm M}_{16}=\lg a,b\mid a^8=b^2=1, bab=a^5\rg.
\end{array}\]
A trivalent symmetric graph $\G$ is said to be of {\em type $2^2$} if $\Aut(\G)_u\cong \S_3$ and $\Aut(\G)_{\{u,v\}}\cong C_{4}$, where $\{u,v\}$ is an edge of $\G$ (see \cite{ConderNedela}). Now we state our next main theorem.

\begin{theorem}\label{mainth2}
Let $\G$ be a locally $2\K_n$ graph with $n\geq 2$. Then $\G$ is $3$-connected-set-homogeneous but not $3$-connected-homogeneous if and only if $\G$ is isomorphic to the line graph $L(\Sigma)$ of a graph $\Sigma$ such that one of the following holds:
\begin{enumerate}
  \item [{\rm (1)}]\  $\Sigma$ is a tetravalent $3$-arc-regular graph;
  \item [{\rm (2)}]\  $\Sigma$ is a pentavalent $3$-arc-regular graph;
  \item [{\rm (3)}]\  $\Sigma$ is a $3$-arc-transitive graph of valency $8$ and $\Aut(\Sigma)_u^{\Sigma(u)}\cong C_2^3\rtimes (C_{7}\rtimes C_{t})$ with $t=1$ or $3$ and $u\in V(\Sigma)$;
  \item [{\rm (4)}]\  $\Sigma$ is a $3$-arc-transitive graph of valency $32$ and $\Aut(\Sigma)_u^{\Sigma(u)}\cong C_2^5\rtimes (C_{31}\rtimes C_{5})$ with $u\in V(\Sigma)$;
  \item [{\rm (5)}]\ $\Sigma$ is a $3$-arc-transitive graph of valency $q+1$ and $\Aut(\Sigma)_u^{\Sigma(u)}\cong\PSL(2,q).\lg\eta\rg$, where $u\in V(\Sigma)$, $q$ is an odd prime power such that $q\equiv-1\ (\mod 4)$ and $\eta$ is a field automorphism of $\GF(q)$;
  \item [{\rm (6)}]\  $\Sigma$ is a pentavalent symmetric graph of type ${\mathcal{Q}_2^6}$;

  \item [{\rm (7)}]\  $\Sigma$ is a trivalent symmetric graph of type $2^2$.
  \end{enumerate}
\end{theorem}

\f{\bf Remark on Theorem~\ref{mainth2}~(5).}\ Let $\Sigma$ be a $3$-arc-transitive graph of valency $q+1$, where $q$ is a prime power. Let $u\in V(\Sigma)$. If $\PGL(2,q)\leq\Aut(\Sigma)_u^{\Sigma(u)}$, then $\Sigma$ is locally $3$-transitive, and then by Theorem~\ref{mainth}, the line graph $\G$ of $\Sigma$ is $3$-connected homogeneous.\smallskip

By Theorem~\ref{mainth2}, to construct $3$-CSH but not $3$-CH graphs which are locally $2\K_n$ with $n\geq 2$, it is equivalent to construct graphs satisfying the conditions in each of (1)--(7) of Theorem~\ref{mainth2}. In \cite[Remark~4.2]{Zhou-EJC-2021}, we gave a pentavalent $3$-arc-regular graph of order $5^3$ (see also Example~\ref{exam:5-val-1}). In 2010, C.H. Li et al. in \cite{LNSS} constructed a tetravalent $3$-arc-regular graph with automorphism group $\PL(2,27)$, and in a recent paper \cite{LiLingMa}, J.J. Li et al. gave another six tetravalent $3$-arc-regular graphs. For trivalent symmetric graphs of type $2^2$, by \cite{Conder-10000} there are only eight such graphs on up to $10000$ vertices, and in 2020, Feng et al.~\cite{Feng-K-M-Yang} constructed an infinite family of trivalent symmetric graphs of type $2^2$. We are not aware of any other $2$-arc-transitive graphs satisfying the conditions in  (1)--(7) of Theorem~\ref{mainth2}.

Our third main theorem provides a useful method to construct $3$-CSH but not $3$-CH graphs which are locally $2\K_n$, and using it, we can give some new constructions of $2$-arc-transitive graphs satisfying the conditions in (1)--(7) of Theorem~\ref{mainth2}. To state the result, we introduce the concept of Cayley graphs.

Given a finite group $G$ and an inverse closed subset $S\subseteq G\setminus\{1\}$, the {\em Cayley graph} $\Cay(G,S)$ on $G$ with respect to $S$ is a graph with vertex set $G$ and edge set $\{\{g,sg\}\mid g\in G,s\in S\}$. For any $g\in G$, $R(g)$ is the permutation of $G$ defined by $R(g): x\mapsto xg$ for $x\in G$.
Set $R(G)=\{R(g)\ |\ g\in G\}$. It is well known that $R(G)$ is a regular subgroup of $\Aut(\Cay(G,S))$. 
In general, a vertex-transitive graph $\G$ is isomorphic to a Cayley graph on a group $G$ if and only if its automorphism group has a subgroup isomorphic to $G$, acting regularly on the vertex set of $\G$ (see \cite[Lemma~16.3]{B}).
Set $A=\Aut(\Cay(G,S))$ and $\Aut(G,S)=\{\a\in\Aut(G)\ |\ S^\a=S\}$. Then $N_A(R(G))=R(G)\rtimes \Aut(G,S)$, and
$\G$ is said to be a {\em normal Cayley graph} of $G$ whenever $N_A(R(G))=A$  (see \cite{Godsil1981,X1}).

Now we state our last theorem which gives a good description for $3$-CSH but not $3$-CH graphs which are locally $2\K_n$ and have solvable automorphism groups. We say that a graph $\G$ is {\em solvable} if $\Aut(\G)$ is solvable.

\begin{theorem}\label{th:sol-3csh-graphs}
Let $n\geq 2$ and let $\G$ be a solvable locally $2\K_n$ graph. Then $\G$ is $3$-connected-set-homogeneous but not $3$-connected-homogeneous if and only if
$\G\cong\Cay(H,S)$ such that the following hold:
\begin{enumerate}
  \item [{\rm (a)}]\ $H$ is a group having two subgroups $A,B$ such that $G=\lg A,B\rg$, $A\cong B\cong C_p^f$, $A\cap B=1$ and $S=(A\cup B)-\{1\}$; and
  \item [{\rm (b)}]\ one of the following holds:
  \begin{enumerate}
    \item [{\rm (1)}]\  $(p,f)=(2,2)$, $\Aut(H,S)\cong C_3\wr C_2$;
   \item [{\rm (2)}]\ $(p,f)=(5,1)$, $\Aut(H,S)\cong C_4\wr C_2$;
   \item [{\rm (3)}]\  $(p,f)=(2,3)$, $C_7\wr C_2\leq \Aut(H,S)\leq (C_7\rtimes C_3)\wr C_2$;
   \item [{\rm (4)}]\ $(p,f)=(2,5)$, $(C_{31}\times C_{31})\rtimes C_{10}\leq \Aut(H,S)\leq (C_{31}\rtimes C_5)\wr C_2$;
  \item [{\rm (5)}]\ $(p,f)=(5,1)$, $\Aut(H,S)\cong {\rm M}_{16}$;
  \item [{\rm (6)}]\ $(p,f)=(3,1)$, $\Aut(H,S)\cong C_4$.
  \end{enumerate}
\end{enumerate}
\end{theorem}

\f{\bf Remark on Theorem~\ref{th:sol-3csh-graphs}}\
{\rm (1)\  Applying Theorem~\ref{th:sol-3csh-graphs}, in Section~\ref{subsec1} we shall show that there are infinitely many graphs $\Cay(H,S)$ satisfying the conditions in each of (1)--(6) of Theorem~\ref{th:sol-3csh-graphs}~(b).

(2)\ By \cite[Theorem~1.1~\&~Corollary~1.2]{LiLingMa}, every tetravalent $3$-arc-regular Cayley graph is a normal cover of a Cayley graph on one of the following groups: $C_3^{11}\rtimes (C_2^{12}.{\rm M}_{11})$, $\S_{35}$ and $\A_{35}$. This, however, is not true. Actually, we shall prove in Proposition~\ref{family-1} there are infinitely many graphs $\G=\Cay(H,S)$ such that $\G$ satisfies the condition in (1) of Theorem~\ref{th:sol-3csh-graphs}~(b), and $\G=L(\Sigma)$ with $\Sigma$ a tetravalent $3$-arc-regular Cayley graph on a solvable group.

(3)\ By \cite{Morgan}, there are nine types of pentavalent $2$-arc-transitive graphs, characterized by the stabilizers of a vertex and an edge. We construct a pentavalent symmetric graphs of type ${\mathcal{Q}_2^6}$. To the best of our knowledge, this is the first known such graph.

(4)\ In 2006, Feng and Kwak~\cite[p.161]{Feng-Kwak} conjectured that every trivalent symmetric graph of order $2\cdot3^n$ is a Cayley graph for each $n\geq 1$. In Lemma~\ref{lem:cubic2-2type}, we shall prove that there exists a Cayley graph $\G=\Cay(H,S)$ on a group $H$ of order $3^{4n+1}$ for each $n\geq 2$ satisfying the condition in (6) of Theorem~\ref{th:sol-3csh-graphs}~(b), and so $\G=L(\Sigma)$, where $\Sigma$ is a trivalent symmetric graph of order $2\cdot 3^{4n}$ of type $2^2$. Then every automorphism of $\Sigma$ swapping any two adjacent vertices of $\Sigma$ is not an involution, and so $\Sigma$ is non-Cayley. This implies that Feng-Kwak's conjecture is not true.}

\medskip

\section{Preliminaries}\label{sec:2}


Let $G$ be a permutation group on a set $\Omega$. For a point $\a\in\Omega$, denote by $G_\a$ the stabilizer of $\a$ in $G$, and denote by $\a^G$ the orbit of $G$ on $\Omega$ containing $\a$. Furthermore, for a subset $\Delta\subseteq \Omega$, denote by $G_\Delta$ the subgroup of $G$ fixing $\Delta$ setwise. If $G$ fixes $\Delta$ setwise, then denote by $G^\Delta$ the permutation group on $\Delta$ induced by $G$.

We say that $G$ is {\em semiregular} on $\Omega$ if $G_\a=1$ for every $\a\in \Omega$ and {\em regular} if $G$ is transitive and semiregular. And, $G$ is said to be {\em primitive} if $G$ is transitive on $\Omega$ and the only partitions of $\Ome$ preserved by $G$ are either the singleton subsets or the whole of $\Ome$. Let $G$ be a transitive permutation group on a set $\Ome$ and let $u\in\Omega$. The orbits of $G_u$ on $\Ome$ are called {\em suborbits} of $G$,
and their sizes are called the {\em subdegrees} of $G$. The number $r$ of the orbits of $G_u$ on $\Ome$ is called the permutation {\em rank} of $G$ on $\Ome$.

A finite transitive permutation group $G$ on a set $\Omega$ is said to be {\em $\frac{3}{2}$-transitive} if all orbits of the stabilizer $G_\a$ of any point $\a\in\Omega$ on $\Omega\setminus\{\a\}$ have the same size greater than $1$. A $\frac{3}{2}$-transitive permutation group $G$ on a set $\Omega$ is said to be a {\em Frobenius group} if $G_{\a\b}=1$ for each different point $\a,\b\in\Omega$. By \cite[Theorems~1.1--1.2]{B-G-L-P-S}, we obtain the following lemma.

\begin{lem}\label{half-transitive}
Let $G$ be a primitive $\frac{3}{2}$-transitive permutation group. Then $G$ is either affine or almost simple. If $G$ is almost simple, then
one of the following holds:
\begin{enumerate}
\item [{\rm (1)}]\ $G$ is $2$-transitive, or
  \item [{\rm (2)}]\ $n=21$, $G=A_7$ or $S_7$ acting on the set of pairs in $\{1,\ldots, 7\}$, or
  \item [{\rm (3)}]\  $n=\frac{1}{2}q(q-1)$ where $q=2^f\geq8$, and either $G=\PSL_2(q)$, or $G={\rm P\G L}_2(q)$ with $f$ prime; the size of
the nontrivial subdegrees is $q + 1$ or $f(q + 1)$, respectively.
\end{enumerate}
\end{lem}

A $2$-arc $(u,v,w)$ of a graph $\G$ is called a {\em $2$-geodesic} if $u$ and $w$ are at distance $2$. A graph $\G$ is said to be {\em $2$-geodesic-transitive} if $\G$ has at least one $2$-geodesic and $\Aut(\G)$ is transitive on the set of $t$-geodesics of $\G$ for $t\leq 2$, where a $1$-geodesic of $\G$ is an arc of $\G$. 

For a $2$-arc $(v_0,v_1,v_2)$ of a graph $\G$, $(v_2,v_1,v_0)$ is also a $2$-arc. If we identify these two arcs, then we obtain a {\em $2$-path}, denoted by $[v_0,v_1,v_2]$, and if $v_0$ and $v_2$ are adjacent then we get a triangle, denoted by $\{v_0,v_1,v_2\}$. The $2$-path $[v_0,v_1,v_2]$ is called a {\em $2$-geodesic-path} provided that the triple $(v_0,v_1,v_2)$ is a $2$-geodesic. We say that $\G$ is {\em $2$-path transitive} ({\em $2$-geodesic-path transitive}, respectively) if $\Aut(\G)$ is transitive on the set of $2$-paths ($2$-geodesic-paths, respectively) of $\G$.

For a graph $\G$, we use $\G^c$ to denote the complementary graph of $\G$.  
From \cite[Theorems~1.1--1.2, Corollary~1.4]{Zhou-EJC-2021}, we obtain the following lemma.

\begin{lem}\label{3-csh-stabilizer}
Let $\G$ be a connected $3$-CSH non-complete graph of girth $3$. Let $G=\Aut(\G)$. Then $\G$ is arc-transitive, and for any $\{u,v\}\in E(\G)$, we have the following:
\begin{enumerate}
\item [{\rm (1)}]\ If $[\G(u)]$ is connected, then $[\G(u)]$ is of diameter $2$, and if $[\G(u)]$ is disconnected, then $[\G(u)]\cong m\K_\ell$ for some positive integers $m,\ell$.
    \item [{\rm (2)}]\ $G_u$ is edge-transitive on $[\G(u)]^c$.
  \item [{\rm (3)}]\ $G_{uv}$ has $s$ orbits on $\G(u)\cap \G(v)$ with equal size, where $s=1, 2, 3$ or $6$.
  \item [{\rm (4)}]\ $G_{uv}$ has $t$ orbits on $\G(u)-((\G(u)\cap \G(v))\cup\{v\})$ with equal size, where $t=1$ or $2$.
  \item [{\rm (5)}]\ If $t=1$, then $\G$ is $2$-geodesic-transitive.
      \end{enumerate}
\end{lem}

The following lemma gives a characterization of $2$-geodesic-path transitive graphs.

\begin{lem}\label{lem:2-geod-path}
Let $\G$ be a connected vertex-transitive graph of valency at least $2$. Take $u\in V(\G)$. Then $\Aut(\G)$ is transitive on the set of $2$-geodesic-paths if and only if $\Aut(\G)_u$ is transitive on the edges of $[\G(u)]^c$.
\end{lem}

\demo Suppose first that $\Aut(\G)$ is transitive on the set of $2$-geodesic-paths. Take two edges, say $\{x,y\}, \{x',y'\}$,  of $[\G(u)]^c$. Then $[x, u, y]$ and $[x',u,y']$ are two $2$-geodesic-paths of $\G$. Then there exists $g\in\Aut(\G)$ sending $[x, u, y]$ to $[x',u,y']$. It follows that $g$ fixes $u$ and sends $\{x,y\}$ to $\{x',y'\}$. Therefore,  $\Aut(\G)_u$ is transitive on the edges of $[\G(u)]^c$.

Conversely, assume that $\Aut(\G)_u$ is transitive on the edges of $[\G(u)]^c$. Take an edge, say $\{x,y\}$, of $[\G(u)]^c$. Then $[x, u, y]$ is  a $2$-geodesic-path of $\G$. For any $2$-geodesic path, say $[x',u',y']$ of $\G$, by the vertex-transitivity, there exists a $g\in\Aut(\G)$ sending $u'$ to $u$, and so $\{x',y'\}^g$ is an edge of $[\G(u)]^c$. Since $\Aut(\G)_u$ is transitive on the edges of $[\G(u)]^c$, there exists $h\in\Aut(\G)_u$ such that $\{x',y'\}^{gh}=\{x,y\}$, and so $[x',u',y']^{gh}=[x, u, y]$. This implies that $\Aut(\G)$ is transitive on the set of $2$-geodesic-paths.\hfill\qed

A {\em clique} of a graph $\G$ is a complete subgraph and a {\em maximal clique} is a clique which is not contained in a larger clique. The {\em clique graph} $C(\G)$ of $\G$ is a graph with vertices the maximal cliques of $\G$ and with two different maximal cliques adjacent if they share at least one common vertex.

\begin{lem}\label{lem:clique-line}
Let $n\geq 2$ be an integer, and let $\G$ be a locally $2\K_n$ graph. Then $\G$ is isomorphic to the line graph of $C(\G)$, and $\Aut(\G)\cong\Aut(C(\G))$.
\end{lem}

\f\demo By {\rm \cite[Corollary~1.6]{DJLP-JCTA}}, we know that $\G$ is isomorphic to the line graph of $C(\G)$, and  by \cite[p.1455]{Babai}, we have $\Aut(\G)\cong\Aut(C(\G))$. \hfill\qed

From {\rm \cite[Theorem~1.1]{DJLP-AMC}} we deduce the following result.

\begin{prop}\label{prop:jinwei-line}
Let $\G$ be a connected regular, non-complete graph of valency at
least $3$. Let $s=2$ or $3$. Then $\G$ is $s$-arc-transitive if and only if the line graph of $\G$ is $(s-1)$-geodesic-transitive.
\end{prop}

\section{Proofs of Theorem~\ref{mainth} and Proposition~\ref{lem:sol-s-tran-stabi}}

In this section, we shall prove Theorem~\ref{mainth} and Proposition~\ref{lem:sol-s-tran-stabi}.\medskip

\f{\bf Proof of Theorem~\ref{mainth}}\ Let $\Sigma=C(\G)$. By Lemma~\ref{lem:clique-line}, for convenience, we shall identify $\G$ with the line graph of $\Sigma$.

Suppose first that $\Sigma$ is $3$-arc-transitive and locally $3$-transitive. By Proposition~\ref{prop:jinwei-line}, $\G$ is $2$-geodesic-transitive and arc-transitive. To show that $\G$ is $3$-CH, it suffices to prove that $\Aut(\G)$ is transitive on the set of $3$-tuples $(e, f, g)$ such that $\{e, f, g\}$ is a triangle. Let $(e_1,e_2,e_3)$ and $(f_1,f_2,f_3)$ be two $3$-tuples of vertices in $\G$ such that both $\{e_1,e_2,e_3\}$ and $\{f_1,f_2,f_3\}$ induce two triangles. As we assume that $\G$ is the line graph of $\Sigma$, we may let $e_1=\{u, v_1\}$, $e_2=\{u,v_2\}$ and $e_3=\{u,v_3\}$, and $f_1=\{x, y_1\}$, $f_2=\{x,y_2\}$ and $f_3=\{x,y_3\}$, where $u,v_i,x,y_i\in V(\Sigma) (i=1,2,3)$. Then $(v_1,u,v_2)$ and $(y_1,x,y_2)$ are $2$-arcs of $\Sigma$. Since $\Sigma$ is $3$-arc-transitive, there exists $\a\in\Aut(\Sigma)$ such that $(v_1,u,v_2)^\a=(y_1,x,y_2)$, and so $(e_1,e_2,e_3)^\a=(f_1,f_2,e_3^\a)$, where $e_3^\a=\{x,v_3^\a\}$.  Clearly, $v_3^\a\in\Sigma(x)$. Since $\Sigma$ is locally $3$-transitive, there exists $\b$ such that $\b$ fixes $x,y_1,y_2$ and maps $v_3^\a$ to $y_3$. So $(e_1,e_2,e_3)^{\a\b}=(f_1,f_2,f_3)$.

Suppose now that $\G$ is $3$-CH. Then $\G$ is $2$-geodesic-transitive, and by Proposition~\ref{prop:jinwei-line}, $\Sigma$ is $3$-arc-transitive. Let $u\in V(\Sigma)$. Then $\Aut(\Sigma)_u$ acts $2$-transitively on $\Sigma(u)$. Take $v,w\in\Sigma(u)$. For any $x_1,x_2\in\Sigma(u)-\{v,w\}$, both $\{\{u,v\},\{u,w\},\{u,x_1\}\}$ and $\{\{u,v\},\{u,w\},\{u,x_2\}\}$ induce two triangles of $\G$. Since $\G$ is $3$-CH, there exists $\a\in\Aut(\Sigma)=\Aut(\G)$ such that $(\{u,v\},\{u,w\},\{u,x_1\})^\a=(\{u,v\},\{u,w\},\{u,x_2\})$. It follows that $\a\in\Aut(\Sigma)_{uvw}$ and $x_1^\a=x_2$. Thus, $\Aut(\Sigma)_u$ acts $3$-transitively on $\Sigma(u)$, and hence $\Sigma$ is locally $3$-transitive. \hfill\qed

\f{\bf Proof of Proposition~\ref{lem:sol-s-tran-stabi}}\  Let $p$ be a prime and $f$ be a positive integer. Suppose that $r$ is a common prime divisor of $p^f-1$ and $f$. Let $\ell$ be a positive integer such that $r^\ell\mid p^f-1$ but $r^{\ell+1}\nmid p^f-1$. Let $\G=\K_{p^f,p^f}$ with biparts $U$ and $W$. Then $\Aut(\G)$ has a subgroup $A=\AGammaL(1,p^f)$ which fixes $U$ point-wise and is $2$-transitive on $W$. Let $g$ be any involution of $\Aut(\G)$ swapping $U$ and $W$, and let $M=\lg A,g\rg$. Then $M=(A\times A^g)\rtimes\lg g\rg$ is $3$-arc-transitive on $\G$, and $M_u^{[1]}$ is transitive on $\G(w)-\{u\}=U-\{u\}$, where $u\in U, w\in W$.

For convenience, we assume that \[A=\AGammaL(1,p^f)=N\rtimes(\lg a\rg\rtimes\lg b\rg)\cong C_p^f\rtimes (C_{p^f-1}\rtimes C_f).\]
Let $B=\lg N, a^{r^\ell}, a^{\frac{p^f-1}{r^\ell}}(b^{\frac{f}{r}})^g,g\rg$. We may let $A_{w}=\lg a\rg\rtimes\lg b\rg$ and $A_{wv}=\lg b\rg$ with $w,v\in W$. Let $u=w^g$ and $x=v^g$. Then $u,x\in U$ and $(w,u,v)$ is a $2$-arc of $\G$. Further, $(A^g)_u=\lg a^g\rg\rtimes\lg b^g\rg$ and $(A^g)_{ux}=\lg b^g\rg$.
Now $B_{wu}=\lg a^{r^\ell}, (a^{r^\ell})^g, a^{\frac{p^f-1}{r^\ell}}(b^{\frac{f}{r}})^g, (a^{\frac{p^f-1}{r^\ell}})^gb^{\frac{f}{r}}\rg$, and
$B_w=N^g\rtimes B_{wu}$ and $B_u=N\rtimes B_{wu}$.  So $B_w$ acts $2$-transitively on $U$. Since $g\in B$, $\G$ is $(B,2)$-arc-transitive. Note that $B_{wuv}=\lg (a^{r^{\ell}})^g, (a^{\frac{p^f-1}{r^\ell}})^gb^{\frac{f}{r}}\rg$ and $B_u^{[1]}=\lg (a^{r^\ell})^g, (a^{\frac{p^f-1}{r^{\ell-1}}})^g\rg$. So $(B_u^{\G(u)})_{wv}\cong B_{wuv}/B_u^{[1]}\leq C_f$. Since $\lg (a^{r^\ell})^g, (a^{\frac{p^f-1}{r^\ell}})^g\rg$ is regular on $U-\{u\}$,  $B_{wuv}$ is also transitive on $\G(v)-\{u\}=U-\{u\}$. It follows that $\G$ is $(B,3)$-arc-transitive. Note that $B_u^{[1]}=\lg (a^{r^\ell})^g, (a^{\frac{p^f-1}{r^{\ell-1}}})^g\rg$ has order $\frac{p^f-1}{r}$. Since $\lg a^g\rg$ is regular on $\G(w)-\{u\}=U-\{u\}$, the number of orbits of $B_u^{[1]}$ on $\G(w)-\{u\}=U-\{u\}$ is $r$, which is a divisor of gcd$(p^f-1,f)$. Now we conclude that $\G$ is a $(G,3)$-arc-transitive graph satisfying the condition in part~(3) of Theorem~\ref{th:3-arc-tran-stabi} with $G=B$.\hfill\qed

\section{Proof of Theorem~\ref{mainth2}}

We begin by proving two lemmas regarding $\frac{3}{2}$-transitive permutation groups.

\begin{lem}\label{impri-half-trans}
Let $G$ be a $\frac{3}{2}$-transitive permutation group on a set $\Omega$. Then the following hold:
\begin{itemize}
  \item [{\rm (1)}]\ if $G$ has rank $3$ or $4$, then $G$ is primitive;
  \item [{\rm (2)}]\ if $G$ has rank $7$, then either $G$ is primitive or $|\Omega|=25$ and $|G|=100$.
\end{itemize}
\end{lem}

\f\demo Suppose that $G$ is imprimitive. By \cite[Theorem~10.4]{WI}, $G$ is a Frobenius group. By \cite[Proposition~4.4]{Sims1967}, there exist two different $x,y\in\Omega$ such that the digraph $\G$ with vertex set $\Omega$ and arc set $\{(x^g,y^g)\mid g\in G\}$ is disconnected.
Let $\G_1$ be a component of $\G$ with $\Omega_1=V(\G_1)$. Then $|\Omega_1|\leq\frac{|\Omega|}{2}$, and since $G$ is transitive on $\Omega$, one has $|\Omega|=t|\Omega_1|$ for some integer $t>0$.

Let $r$ be the rank of $G$ and take $\a\in\Omega_1$. Since $G$ is a Frobenius group, one has $|G_\a|=\frac{|\Omega|-1}{r-1}\geq 2$. It implies that $|\Omega|\geq 2r-1$.
Clearly, $G_\a$ fixes $\Omega_1$ setwise. Again, since $G$ is a Frobenius group, we know that $G_\a$ acts semiregularly on $\Omega_1-\{\a\}$, and so $|\Omega_1|=1+k|G_\a|$ for some positive integer $k<r-1$. It follows that $|\Omega_1|=1+\frac{k(|\Omega|-1)}{r-1}.$
Since $|\Omega|=t|\Omega_1|$, one has $|\Omega|=t\cdot(1+\frac{k(|\Omega|-1)}{r-1})$, implying that
\begin{equation}\label{eq-3}
(r-1-tk)|\Omega|=t(r-1-k)>0.
\end{equation}
It follows that $r-1-tk>0$. If $t=1$, then Eq.~(\ref{eq-3}) implies that $|\Omega|=1$, contrary to $|\Omega|\geq 2r-1$. Thus, $t>1$ and hence $r>2k+1$.

This implies that $r>3$. If $r=4$, then $k=1$, and then by Eq.~(\ref{eq-3}), we obtain that $(3-t)|\Omega|=2t.$ It implies that $t<3$ and $|\Omega|\leq 4$, contrary to $|\Omega|\geq 2r-1=7$. This proves part (1).



If $r=7$, then $2k+1<7$, and then $k\leq 2$. In case $k=2$,  by Eq.~(\ref{eq-3}), we have $(6-2t)|\Omega|=4t.$ It follows that $t<3$ and $|\Omega|\leq 4$, contrary to $|\Omega|\geq 2r-1=13$. Thus, we have $k=1$. Again, by Eq.~(\ref{eq-3}), we obtain that $(6-t)|\Omega|=5t$, implying $t\leq 5$. Since $6=r-1$ is a divisor of $|\Omega|-1$, it follows that $t=5$, and hence $|\Omega|=25$. Then $|G|=100$. Part (3) holds.\hfill\qed

\begin{lem}\label{lem:2-tran+impri}
Let $G$ be a $2$-transitive permutation group on a set $\Omega$ and take $x\in\Omega$. Suppose that $G_x$ is a $\frac{3}{2}$-transitive permutation group on $\Omega-\{x\}$ of rank $3$, $4$ or $7$. Then $G_x$ is primitive on $\Omega-\{x\}$ and one of the following holds.
\begin{enumerate}
  \item [{\rm (1)}]\ $G=\PSL(2,q).\lg\eta\rg$ and $|\Omega|=q+1$, where $q>3$ is an odd prime power and $\eta$ is a field automorphism of $\GF(q)$. Moreover, $G_x\cong C_{\frac{q-1}{2}}.\mathcal{O}$ has rank $3$ on $\Omega-\{x\}$, where $\mathcal{O}\cong\lg\eta\rg$.
  \item [{\rm (2)}]\ $G\cong C_2^f\rtimes(C_{2^f-1}\rtimes C_f)$ and $|\Omega|=2^f$, where $f=3$ or $5$. Moreover, $G_x$ has rank $3$ or $7$ on $\Omega-\{x\}$.
\end{enumerate}
\end{lem}

\f\demo Suppose first that $G_x$ is imprimitive on $\Omega-\{x\}$. By Lemma~\ref{impri-half-trans}, $G_x$ is a permutation group on $\Omega-\{x\}$ with rank $7$, and $|\Omega-\{x\}|=25$ and $|G_x|=100$. However, by checking \cite[Appendix~B]{Dixon-Mortimer}, there exist no $2$-transitive permutation groups of degree $26$ with point-stabilizer of order $100$, a contradiction. Thus, $G_x$ is primitive on $\Omega-\{x\}$. 

Since $G_x$ is $\frac{3}{2}$-transitive on $\Omega-\{x\}$ with rank $3$, $4$ or $7$, it follows that $G$ is $\frac{5}{2}$-transitive but neither $3$-transitive nor sharply $2$-transitive on $\Omega$. By \cite[Proposition~4]{Libeck-Praeger-Saxl2019}, we conclude that $G$ is one of the groups in Table~\ref{tab-3}.
\begin{table}
\begin{center}
\caption{2-Transitive groups}\label{tab-3}
\begin{tabular}{lccl}
\hline Line & $G$ &   $|\Omega|$ &  Remarks \\
\hline
1 & $\PSL(2,q)\unlhd G\leq {\rm P\Gamma L}(2, q)$ & $q+1$ & $q>3$ a prime power\\

2 & {\rm Sz}$(q)$  & $q^2+1$ & $q=2^{2a+1}>2$\\

3 & ${\rm A\Gamma L}(1,2^f)$  & $2^f$ & $f$ a prime\\
\hline

\end{tabular}
\end{center}
\end{table}

If $\PSL(2, q)\unlhd G\leq {\rm P\Gamma L}(2,q)$, then $|\Omega|=q+1$ and $[q]: C_{\frac{q-1}{s}}\unlhd G_x$ with $s=(2, q-1)$.
Since $G$ is not $3$-transitive on $\Omega$, one has $\PGL(2, q)\nleq G$ and $q>3$ is odd. It follows that $G=\PSL(2,q).\lg\eta\rg$, where $\eta$ is a field automorphism of $\GF(q)$. Moreover, $G_x\cong C_{\frac{q-1}{2}}.\mathcal{O}$ with $\mathcal{O}\cong\lg\eta\rg$, and $G_x$ has rank $3$ on $\Omega-\{x\}$. So part (1) holds.

If $G={\rm Sz}(q)$, then $G_x\cong [q^2]:C_{q-1}$. By \cite[Section~7.7]{Dixon-Mortimer}, we deduce that the normal subgroup of $G_x$ of order $q^2$ is not an elementary abelian $2$-group. This contradicts that $G_x$ is primitive on $\Omega-\{x\}$.

Now let $G={\rm A\Gamma L}(1, 2^f)$ with $f$ a prime. Then $|\Omega|=2^f$ and $G_x\cong C_{2^f-1}\rtimes C_f$. Since $G_x$ is primitive on $\Omega-\{x\}$, one has $|\Omega-\{x\}|=2^f-1$ and $\soc(G_x)\cong C_{2^f-1}$. It follows that $2^f-1$ is a prime.
Furthermore, for an arbitrary $y\in \Omega-\{x\}$, we have $G_{xy}\cong C_f$ and $G_{xy}$ is semiregular on $\Omega-\{x,y\}$. Since $G_x$ is $\frac{3}{2}$-transitive on $\Omega-\{x\}$ with rank $3$, $4$ or $7$, it follows that $G_{xy}$ has $s$ orbits of size $f$ on $\Omega-\{x,y\}$, where $s=2,3$ or $6$. This implies that $sf=|\Omega-\{x,y\}|=2^f-2$.

Recall that $f$ is a prime. If $f>5$, then $s=2$ or $6$, and hence $2^{f-1}-1=f$ or $3f$. However, this is impossible because $f>5$. If $f=2$, then $G_x\cong C_3\rtimes C_2$ which is $2$-transitive on $\Omega-\{x\}$. This is impossible because $G_x$ is $\frac{3}{2}$-transitive group on $\Omega-\{x\}$. Thus, $f=3$ or $5$, and $G_x\cong C_{2^f-1}\rtimes C_f$. Part (2) happens.\hfill\qed

In the next four lemmas, we shall prove the sufficiency of Theorem~\ref{mainth2}.

\begin{lem}\label{lem:3CSH-2arc}
Let $\G$ be a connected $3$-CSH locally $2\K_n$ graph with $n\geq 2$. Then $C(\G)$ is $2$-arc-transitive.
\end{lem}

\f\demo By Lemma~\ref{lem:clique-line}, $\G$ is isomorphic to the line graph of $C(\G)$, and by Lemma~\ref{3-csh-stabilizer}, $\G$ is arc-transitive and so $1$-geodesic-transitive. It then follows from Proposition~\ref{prop:jinwei-line} that $C(\G)$ is $2$-arc-transitive.\hfill\qed

\begin{lem}\label{lem:triangle-tran-1}
Let $\G$ be a connected locally $2\K_{q}$ graph, where $q>3$ is an odd prime power. Let $\Sigma=C(\G)$ and let $A=\Aut(\Sigma)$. Take $u\in V(\Sigma)$ and $v\in\Sigma(u)$. Suppose that $\Sigma$ is $3$-arc-transitive and that $A_u^{\Sigma(u)}\cong\PSL(2,q).\lg\eta\rg$, where $\eta$ is a field automorphism of $\GF(q)$. Then  $\G$ is $3$-CSH if and only if $q\equiv-1\ (\mod 4)$.
\end{lem}

\f\demo By Lemma~\ref{lem:clique-line}, we may assume that $\G$ is the line graph of $\Sigma$. Take an edge $e=\{u_e, v_e\}$ of $\Sigma$. The number of triangles of $\G$ passing through $e$ is equal to $q(q-1)$, which is just the number of edges of $[\G(e)]$.  Since $\Sigma$ is $3$-arc-transitive, $\G$ is vertex-transitive, and as every triangle contains three vertices, the total number of triangles of $\G$ is $q(q-1)|V(\G)|/3$.

\begin{figure}[ht]
\begin{center}
\unitlength 4mm
\begin{picture}(10,7)

\put(3, 7){\circle*{0.4}} \put(1.7, 7){$e_1$}
\put(3, 5){\circle*{0.4}} \put(1.7, 5){$e_2$}
\put(3, 3){$\vdots$} 
\put(3, 1){\circle*{0.4}} \put(1.7, 1){$e_q$}

\put(3,7){\line(1, 0){4}}
\put(3,7){\line(2, -1){4}}
\put(3,7){\line(2, -3){4}}
\put(3,5){\line(1, 0){4}}
\put(3,5){\line(2, 1){4}}
\put(3,5){\line(1, -1){4}}
\put(3,1){\line(2, 3){4}}
\put(3,1){\line(1, 1){4}}
\put(3,1){\line(1, 0){4}}

\put(7, 7){\circle*{0.4}} \put(7.7, 7){$e_1'$}
\put(7, 5){\circle*{0.4}} \put(7.7, 5){$e_2'$}
\put(7, 3){$\vdots$} 
\put(7, 1){\circle*{0.4}} \put(7.7, 1){$e_q'$}
\end{picture}
\end{center}\vspace{-.5cm}
\caption{The complement $[\G(e)]^c$ of $[\G(e)]$} \label{fig-3}
\end{figure}
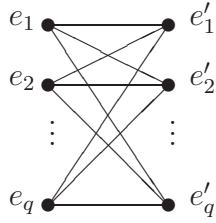

Let $\G(e)=\{e_i, e'_i\mid 1\leq i\leq q\}$, where $e_1, e_2, \cdots, e_q$ are edges of $\Sigma$ incident with $u_e$ and $e'_1, e'_2, \cdots, e'_q$ are edges of $\Sigma$ incident with $v_e$. Then $[\G(e)]\cong 2\K_q$, and $\{e_i\mid i=1,2,\dots,q\}$ and $\{e_i'\mid i=1,2,\dots,q\}$ are the bi-parts of $[\G(e)]^c=\K_{q,q}$ (see Figure~(\ref{fig-3})).

Let $N_1=\{e, e_i\ |\ 1\leq i\leq q\}$. Then $N_1$ is just the set of edges of $\Sigma$ incident with $u_e$ and $|N_1|=q+1$. So we may view $A_{u_e}^{\Sigma(u_e)}$ as a permutation group on $N_1=\{e,e_i\ |\ 1\leq i\leq q\}$.

Note that $|\PSL(2,q)|=\frac{1}{2}q(q^2-1)$. This implies that $3\mid |\PSL(2,q)|$. Take an element, say $x$, of order $3$ in $\soc(A_{u_e}^{\Sigma(u_e)})$. We may assume that $e_1^x\neq e_1$. Then $\{e_1,e_1^x,e_1^{x^2}\}$ induces a triangle of $\G$, and so $C_3\leq A_{\{e_1, e_1^x, e_1^{x^2}\}}/A_{e_1e_1^xe_1^{x^2}}\leq \S_{3}$.
Since $\Sigma$ is $3$-arc-transitive, we obtain that $\G$ is arc-transitive. It follows that $|A:A_{e_1}|=|V(\G)|$ and $|A_{e_1}:A_{e_1e_1^x}|=|\G(e_1)|=2q$.
Noticing that $A_{e_1e_1^x}$ fixes $e_1\cap e_1^x=\{u_e\}$, we have $A_{e_1e_1^x}\leq A_{u_e}$. It implies that $A_{e_1e_1^x}$ fixes $N_1$ setwise as $N_1$ is the set of edges of $\Sigma$ incident with $u_e$. So $|A_{e_1e_1^x}: A_{e_1e_1^xe_1^{x^2}}|=|(A_{u_e}^{\Sigma(u_e)})_{e_1e_1^x}: (A_{u_e}^{\Sigma(u_e)})_{e_1e_1^xe_1^{x^2}}|=\frac{q-1}{2}$.
It follows that \[|A: A_{e_1e_1^xe_1^{x^2}}|=|A:A_{e_1}||A_{e_1}:A_{e_1e_1^x}||A_{e_1e_1^x}: A_{e_1e_1^xe_1^{x^2}}|=q(q-1)|V(\G)|.\]

If $q\equiv 1\ (\mod 4)$, then there exists an involution $y\in (\soc(A_{u_e}^{\Sigma(u_e)}))_{e_1e_1^x}$ and so $y$ must interchange another two edges, say $e_i$ and $e_j$ in $N_1$. This implies that $C_2\leq A_{\{e_1, e_i, e_j\}}/A_{e_1e_ie_j}$. If $\G$ is $3$-CSH, then $A$ is transitive on the triangles of $\G$, and then we would have $A_{\{e_1, e_1^x, e_1^{x^2}\}}/A_{e_1e_1^xe_1^{x^2}}\cong \S_{3}$. It then follows that $|A_{\{e_1, e_1^x, e_1^{x^2}\}}|=6|A_{e_1e_1^xe_1^{x^2}}|$. Consequently, the size of the orbit $\{e_1, e_1^x, e_1^{x^2}\}^A$ of $A$ acting on the set of triangles of $\G$ is
\[|A: A_{\{e_1, e_1^x, e_1^{x^2}\}}|=\frac{q(q-1)}{6}|V(\G)|.\]
This, however, is impossible because the total number of triangles of $\G$ is $q(q-1)|V(\G)|/3$.

If $q\equiv -1\ (\mod 4)$, then $(A_{u_e}^{\Sigma(u_e)})_{e_1e_1^x}$ has odd order, and so every involution in $A_{u_e}^{\Sigma(u_e)}$ does not fix any edge in $N_1$. Since $A_{\{e_1, e_1^x, e_1^{x^2}\}}$ fixes $e_1\cap e_1^x\cap e_1^{x^2}=\{u_e\}$, we have $A_{e_1e_1^x}\leq A_{u_e}$. It follows that $A_{\{e_1, e_1^x, e_1^{x^2}\}}$ fixes $N_1$ setwise since $N_1$ is the set of edges of $\Sigma$ incident with $u_e$. So $A_{\{e_1, e_1^x, e_1^{x^2}\}}^{\G(e)}\leq A_{u_e}^{\Sigma(u_e)}$. This implies that $A_{\{e_1, e_1^x, e_1^{x^2}\}}/A_{e_1e_1^xe_1^{x^2}}\cong C_{3}$. Consequently, the size of the orbit $\{e_1, e_1^x, e_1^{x^2}\}^A$ of $A$ acting on the set of triangles of $\G$ is
\[|A: A_{\{e_1, e_1^x, e_1^{x^2}\}}|=\frac{q(q-1)}{3}|V(\G)|.\]
Thus, $\G$ is $3$-CSH.\hfill\qed

\begin{lem}\label{lem:triangle-tran}
Let $\G$ be a connected locally $2\K_{q}$ graph, where either $q>3$ is an odd prime power or $q=2^f-1$ with $f=2,3$ or $5$. Let $\Sigma=C(\G)$ and let $A=\Aut(\Sigma)$. Take $u\in V(\Sigma)$ and $v\in\Sigma(u)$. Suppose that $\Sigma$ is $3$-arc-transitive and satisfies one of the following:
\begin{enumerate}
   \item [{\rm (1)}]\ $q$ is an odd prime power such that $q\equiv-1\ (\mod 4)$, and $A_u^{\Sigma(u)}\cong\PSL(2,q).\lg\eta\rg$, where $\eta$ is a field automorphism of $\GF(q)$, or,
  \item [{\rm (2)}]\  $q=2^2-1$ and $A_u^{\Sigma(u)}\cong\A_4$, or,
  \item [{\rm (3)}]\ $q=2^3-1$ and $A_u^{\Sigma(u)}\cong C_2^3\rtimes (C_{7}\rtimes C_{s})$ with $s=1$ or $3$, or,
  \item [{\rm (4)}]\ $q=2^5-1$ and $A_u^{\Sigma(u)}\cong C_2^5\rtimes (C_{31}\rtimes C_{5})$.
\end{enumerate}
Then $\G$ is $2$-geodesic-transitive and $3$-CSH but not $3$-CH.
\end{lem}

\f\demo In view of Lemma~\ref{lem:clique-line}, for convenience, we shall assume that $\G$ is the line graph of $\Sigma$. Since $\Sigma$ is $3$-arc-transitive, by Proposition~\ref{prop:jinwei-line}, $\G$ is $2$-geodesic-transitive. Observe that in each of (1)--(4), $A_u^{\Sigma(u)}$ is not $3$-transitive. It implies that $\Sigma$ is not locally $3$-transitive. So $\Gamma$ is not $3$-CH by Theorem~\ref{mainth}.

If $\Sigma$ satisfies the condition in part (1), then by Lemma~\ref{lem:triangle-tran-1}, $\G$ is $3$-CSH.

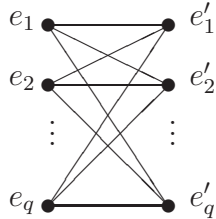
\begin{figure}[ht]
\begin{center}
\unitlength 4mm
\begin{picture}(10,7)

\put(3, 7){\circle*{0.4}} \put(1.7, 7){$e_1$}
\put(3, 5){\circle*{0.4}} \put(1.7, 5){$e_2$}
\put(3, 3){$\vdots$} 
\put(3, 1){\circle*{0.4}} \put(1.7, 1){$e_q$}

\put(3,7){\line(1, 0){4}}
\put(3,7){\line(2, -1){4}}
\put(3,7){\line(2, -3){4}}
\put(3,5){\line(1, 0){4}}
\put(3,5){\line(2, 1){4}}
\put(3,5){\line(1, -1){4}}
\put(3,1){\line(2, 3){4}}
\put(3,1){\line(1, 1){4}}
\put(3,1){\line(1, 0){4}}

\put(7, 7){\circle*{0.4}} \put(7.7, 7){$e_1'$}
\put(7, 5){\circle*{0.4}} \put(7.7, 5){$e_2'$}
\put(7, 3){$\vdots$} 
\put(7, 1){\circle*{0.4}} \put(7.7, 1){$e_q'$}
\end{picture}
\end{center}\vspace{-.5cm}
\caption{The complement $[\G(e)]^c$ of $[\G(e)]$} \label{fig-2}
\end{figure}

Next we consider the case when $\Sigma$ satisfies the condition in one of parts (2), (3) and (4). Let $q=2^f-1$ with $f=2,3$ or $5$. Let $e=\{u,v\}$ and let $\G(e)=\{e_i,e_i'\mid i=1,2,\dots,q\}$, where $e_1, e_2, \cdots, e_q$ are edges of $\Sigma$ incident with $u$ and $e'_1, e'_2, \cdots, e'_q$ are edges of $\Sigma$ incident with $v$. Then $[\G(e)]\cong 2\K_q$, and $\{e_i\mid i=1,2,\dots,q\}$ and $\{e_i'\mid i=1,2,\dots,q\}$ are the bi-parts of $[\G(e)]^c=\K_{q,q}$ (see Figure~(\ref{fig-2})).

Assume first that $\Sigma$ satisfies the condition in part (2). Then $q=3$ and $A_u^{\Sigma(u)}\cong\A_4$. By \cite[Theorem~4]{Potocnik}, we have that $A_u\cong A_4\times C_3$. Then the arc stabilizer $A_{(u, v)}\cong C_3\times C_3$, and $A_{(u, v)}^{\G(e)}=\lg (e_1, e_2, e_3)\rg\times\lg(e_1', e_2', e_3')\rg$. Since $\Sigma$ is arc-transitive, there exists an involution $g\in A_{e}$ such that $g$ swaps the two bi-parts of $[\G(e)]^c=\K_{3,3}$.
This implies that $A_e\cong(C_3\times C_3)\rtimes C_2$ and that $A_e$ is edge-transitive on both $[\G(e)]$ and $[\G(e)]^c$. By \cite[Lemma~3.1]{Zhou-EJC-2021}, $\G$ is $3$-CSH.

Now assume that $\Sigma$ satisfies the condition in one of parts (3) and (4). Then $q=2^f-1$ with $f=3$ or $5$, and $A_{u}^{\Sigma(u)}\cong C_2^f\rtimes (C_{q}\rtimes C_{s})$, where $s=1$ or $3$ if $f=3$, and $s=5$ if $f=5$. To complete the proof, it suffices to show that $A$ is transitive on the triangles of $\G$. Note that the number of triangles of $\G$ passing through $e$ is equal to $q(q-1)$, which is just the number of edges of $[\G(e)]$. Since $\G$ is vertex-transitive and every triangle contains three vertices, the total number of triangles of $\G$ is $q(q-1)|V(\G)|/3$.

Let $N_1=\{e, e_i\ |\ 1\leq i\leq q\}$. Then $N_1$ is just the set of edges of $\Sigma$ incident with $u$ and $|N_1|=q+1$. So we may view $A_{u}^{\Sigma(u)}$ as a permutation group on $N_1$. Then $A_u^{[1]}$ is the kernel of $A_u$ acting on $N_1$ and $A_u^{\Sigma(u)}\cong A_u/A_u^{[1]}$. For any distinct $e_i, e_j\in N_1$, $\{e, e_i, e_j\}$ induces a triangle of $\G$. Since $\Sigma$ is $3$-arc-transitive, it follows that $A_u^{\Sigma(u)}\cong C_2^f\rtimes (C_{q}\rtimes C_{s})$ is $2$-transitive on $\Sigma(u)$. So $A_u^{\Sigma(u)}$ is also $2$-transitive on $N_1$. It follows that $(A_{u}^{\Sigma(u)})_\a\cong C_q\rtimes C_s$ for all $\a\in\{e, e_i, e_j\}$, $(A_{u}^{\Sigma(u)})_{ee_i}\cong C_s$ and $(A_{u}^{\Sigma(u)})_{ee_ie_j}=1$.

Noting that $e\cap e_i\cap e_j=\{u\}$, we have $A_{\{e, e_i, e_j\}}\leq A_{u}$. This implies that $A_{\{e, e_i, e_j\}}$ fixes $N_1$ setwise. As $(A_{u}^{\Sigma(u)})_{ee_ie_j}=1$, one has $A_{ee_ie_j}=A_u^{[1]}\cap A_{\{e, e_i, e_j\}}$. It follows that $A_{\{e, e_i, e_j\}}/A_{ee_ie_j}\cong A_{\{e, e_i, e_j\}}A_u^{[1]}/A_{u}^{[1]}\leq A_{u}^{\Sigma(u)}$. As $(A_{u}^{\Sigma(u)})_\a\cong C_q\rtimes C_s$ for all $\a\in\{e, e_i, e_j\}$, it follows that $2\nmid |A_{\{e, e_i, e_j\}}/A_{ee_ie_j}|$ since $qs$ is odd. Thus, $A_{\{e, e_i, e_j\}}/A_{ee_ie_j}\leq C_3$.

If $(s,f)\neq (3,3)$, then $3\nmid |A_{u}^{\Sigma(u)}|$ and so $3\nmid |A_{\{e, e_i, e_j\}}/A_{ee_ie_j}|$. Then $|A_{\{e, e_i, e_j\}}/A_{ee_ie_j}|=1$.

If $(s,f)=(3,3)$, then we may take an element $x$ of $A_{u}^{\Sigma(u)}$ of order $3$ such that $e^x\neq e$. By the arbitrariness of $e_i$ and $e_j$, we may assume that $\{e, e^x, e^{x^2}\}=\{e, e_i, e_j\}$. Then we have $A_{\{e, e_i, e_j\}}/A_{ee_ie_j}\cong C_3$.

Now $|A: A_{\{e, e_i, e_j\}}|=\frac{1}{3}|A: A_{ee_ie_j}|$ when $s=f=3$, and otherwise, $|A: A_{\{e, e_i, e_j\}}|=|A: A_{ee_ie_j}|$.
As $e\cap e_i=\{u\}$, we have $A_{ee_i}\leq A_{u}$, and so $|A_{ee_i}: A_{ee_ie_j}|=|(A_{u}^{\Sigma(u)})_{ee_i}: (A_{u}^{\Sigma(u)})_{ee_ie_j}|=s$. Since $\G$ is arc-transitive, one has $|V(\G)|=|A:A_e|$ and $2q=|\G(e)|=|A_e: A_{ee_i}|$. It follows that
\[|A: A_{ee_ie_j}|=|A:A_{e}||A_{e}:A_{ee_i}||A_{ee_i}: A_{ee_ie_j}|=2qs|V(\G)|.\]
As a result, the size of the orbit $\{e,e_i,e_j\}^A$ of $A$ acting on the set of triangles of $\G$ is
\[|A: A_{\{e, e_i, e_j\}}|=\frac{q(q-1)}{3}|V(\G)|.\]
It follows that $A$ is transitive on the set of triangles of $\G$.
\hfill\qed

\begin{lem}\label{lem:3-ar-2gt}
Let $\G$ be a connected locally $2\K_n$ graph with $n\geq 2$. If $C(\G)$ is a $3$-arc-regular graph of valency $5$, then $\G$ is $2$-geodesic-transitive and $3$-CSH but not $3$-CH.
\end{lem}

\f\demo Let $\Sigma=C(\G)$. In view of Lemma~\ref{lem:clique-line}, for convenience, we shall assume that $\G$ is the line graph of $\Sigma$. Let $A=\Aut(\Sigma)$. By Proposition~\ref{prop:jinwei-line}, $\G$ is $2$-geodesic-transitive. Take $u\in V(\Sigma)$ and $v\in\Sigma(u)$.
Since $\Sigma$ is a $3$-arc-regular graph of valency $5$, by \cite[Theorem~4.1]{Zhou-Feng-2010}, we have $A_u\cong (C_5\rtimes C_4)\times C_4$ and $A_{\{u,v\}}\cong (C_4\times C_4)\rtimes C_2$. By \cite[Theorem~1.5(2)]{Zhou-EJC-2021}, $\G$ is $3$-CSH not $3$-CH.\hfill\qed

\begin{lem}\label{lem:not2gt}
Let $\G$ be a connected locally $2\K_n$ graph with $n\geq 2$. If $C(\G)$ is either a trivalent symmetric graph of type $2^2$ or a pentavalent symmetric graph of type $\mathcal{Q}_2^6$, then $\G$ is $3$-CSH but not $2$-geodesic-transitive.
\end{lem}

\f\demo Let $\Sigma=C(\G)$. In view of Lemma~\ref{lem:clique-line}, for convenience, we shall assume that $\G$ is the line graph of $\Sigma$.

Assume first that $\Sigma$ is a trivalent symmetric graph of type $2^2$. Then $\Sigma$ is $2$-arc-regular with edge-stabilizer isomorphic to $C_4$. By \cite[Theorem~1.5(2)]{Zhou-EJC-2021}, $\G$ is $3$-CSH, and by Proposition~\ref{prop:jinwei-line}, $\G$ is not $2$-geodesic-transitive.

Assume now that $\Sigma$ is a pentavalent symmetric graph of type $\mathcal{Q}_2^6$. Let $A=\Aut(\Sigma)$. Take an edge $e=\{u_e, v_e\}$ of $\Sigma$. Then $\Sigma$ is $2$-arc-transitive with vertex stabilizer $A_{u_e}\cong\F(20)\times C_2$ and edge stabilizer $A_e\cong {\rm M}_{16}$ (see \cite[Table~1]{Morgan}),
where\[{\rm M}_{16}=\lg a,b\mid a^8=b^2=1, bab=a^5\rg.\]
Again, by Proposition~\ref{prop:jinwei-line}, $\G$ is arc-transitive but not $2$-geodesic-transitive. To complete the proof, it suffices to prove that $\G$ is both $2$-geodesic-path transitive and triangle transitive.

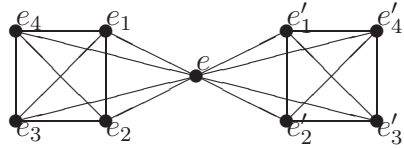
\begin{figure}[ht]
\begin{center}
\unitlength 4mm
\begin{picture}(12,5)

\put(6, 3){\circle*{0.4}} \put(6, 3.2){$e$}

\put(6,3){\line(2, 1){3}} \put(6,3){\line(2, -1){3}}
\put(6,3){\line(-2, 1){3}} \put(6,3){\line(-2, -1){3}}

\put(6,3){\line(4, 1){6}} \put(6,3){\line(4, -1){6}}
\put(6,3){\line(-4, 1){6}} \put(6,3){\line(-4, -1){6}}

\put(3, 4.5){\circle*{0.4}} \put(3, 4.7){$e_1$}
\put(3, 1.5){\circle*{0.4}} \put(3, 1){$e_2$}
\put(0, 4.5){\circle*{0.4}} \put(0, 4.7){$e_4$}
\put(0, 1.5){\circle*{0.4}} \put(0, 1){$e_3$}

\put(0,1.5){\line(0, 1){3}} \put(0,1.5){\line(1, 1){3}}
\put(0,1.5){\line(1, 0){3}}
\put(3,4.5){\line(0, -1){3}} \put(3,1.5){\line(-1, 1){3}}
\put(3,4.5){\line(-1, 0){3}}

\put(9, 4.5){\circle*{0.4}} \put(9, 4.7){$e_1'$}
\put(9, 1.5){\circle*{0.4}} \put(9, 1){$e_2'$}
\put(12, 4.5){\circle*{0.4}} \put(12, 4.7){$e_4'$}
\put(12, 1.5){\circle*{0.4}} \put(12, 1){$e_3'$}

\put(9,1.5){\line(0, 1){3}} \put(9,1.5){\line(1, 1){3}}
\put(9,1.5){\line(1, 0){3}}
\put(12,4.5){\line(0, -1){3}} \put(12,1.5){\line(-1, 1){3}}
\put(12,4.5){\line(-1, 0){3}}
\end{picture}
\end{center}\vspace{-.5cm}
\caption{The subgraph of $\G$ induced by $\{e\}\cup\G(e)$} \label{fig-1}
\end{figure}

Let $\G(e)=\{e_i, e'_i\ |\ 1\leq i\leq 4\}$, where $e_1, e_2, e_3, e_4$ are the four edges of $\Sigma$ incident with $u_e$ and $e'_1, e'_2, e'_3, e'_4$ are the other four edges of $\Sigma$ incident with $v_e$.
Then $[\G(e)]\cong 2\K_4$, see Figure~\ref{fig-1}.
Since $A_{u_e}\cong\F(20)\times C_2$ and $\Sigma$ is $2$-arc-transitive of valency $5$, we have $A_{u_e}^{[1]}\cong C_2$ and $A_{u_e}^{[1]}\cap A_{v_e}^{[1]}=1$. It follows that $A_e$ acts faithfully on $\G(e)$. Since $A_e\cong {\rm M}_{16}$, without loss of generality, we may assume that $A_e={\rm M}_{16}.$ Then $\lg a\rg$ acts regularly on $\G(e)$. Let $B_0=\{e_i\mid 1\leq i\leq 4\}$ and $B_1=\{e_i'\mid 1\leq i\leq 4\}$. Then $B_0$ and $B_1$ are the two bi-parts of $[G(e)]^c\cong\K_{4,4}$, where $[\G(e)]^c$ is the complement of the induced subgraph $[\G(e)]$ of $\G$. The subgroup of $\lg a\rg$ fixing $B_0$ setwise is $\lg a^2\rg$. So we may view $[G(e)]^c$ as the Cayley graph $\Delta=\Cay(\lg a\rg, \{a, a^3, a^5, a^7\})$ on $\lg a\rg$. Furthermore, we may identify $B_0$ with $\lg a^2\rg$ and identify $B_1$ with $a\lg a^2\rg$.

Then $\lg a\rg$ acts on $V(\Delta)=\lg a\rg$ by right multiplication and $\lg b\rg$ acts on $V(\Delta)=\lg a\rg$ by conjugation. Then \[E(\Delta)=\{1, a\}^{\lg a\rg}\cup \{1, a^3\}^{\lg a\rg}\cup\{1, a^5\}^{\lg a\rg}\cup \{1, a^7\}^{\lg a\rg}.\]
Note that $\{1,a\}^{a^{7}}=\{1,a^7\},$ $\{1,a^3\}^{a^{5}}=\{1,a^5\}$ and $\{1,a\}^{b}=\{1,a^5\}$. It follows that $E(\Delta)=\{1, a\}^{A_e}$. This implies that $A_e={\rm M}_{16}$ is transitive on the edges of $[\G(e)]^c$. By Lemma~\ref{lem:2-geod-path}, $\G$ is $2$-geodesic-path transitive.

Note that the subgroup of $A_e$ fixing $B_0$ setwise is $\lg a^2\rg\times\lg b\rg$, which induces a regular action on $B_0$. It follows that for each $1\leq i\leq 4$, $A_{ee_i}$ fixes all $e_1, e_2, e_3, e_4$. So $A_{ee_1}=A_{ee_1e_2}$. 
Consider the triangle of $\G$ induced by $\{e, e_1, e_2\}$. If there exists $g\in A$ such that $g$ cyclically permutes $e, e_1$ and $e_2$, then $g$ must fix the maximal clique of $\G$ induced by $e_1, e_2, e_3, e_4, e$. As $u_e$ is the intersection of these five edges, we have $g\in A_{u_e}$, forcing that $A_{u^e}^{\Sigma(u_e)}$ would contain an element of order $3$. This, however, is impossible because $A_{u_e}\cong \F(20)\times C_2$. Thus, $A_{\{e, e_1, e_2\}}/A_{ee_1e_2}\cong C_k$ with $k=1$ or $2$. Since $A_{ee_1}=A_{ee_1e_2}$, one has $|A_{e_1e_2}: A_{e_1e_2e_3}|=1$. Since $\G$ is arc-transitive, one has $|A:A_e|=|V(\G)|$ and $|A_{e}:A_{ee_1}|=|\G(e)|=8$. It follows that \[|A: A_{\{e_1, e_2, e_3\}}|=\frac{1}{k}|A:A_{e_1}||A_{e_1}:A_{e_1e_2}||A_{e_1e_2}: A_{e_1e_2e_3}|=\frac{8}{k}|V(\G)|\geq 4|V(\G)|.\] The number of triangles of $\G$ is $24|V(\G)|/6=4|V(\G)|$. As $|\{e_1, e_2, e_3\}^A|=|A: A_{\{e_1, e_2, e_3\}}|\leq 4|V(\G)|$, it follows that $|\{e_1, e_2, e_3\}^A|=4|V(\G)|$ and so $\G$ is triangle transitive. This completes the proof. \hfill\qed

Now we are ready to prove Theorem~\ref{mainth2}. Due to Lemma~\ref{lem:clique-line}, it is enough to prove the following theorem.

\begin{theorem}\label{mainth:locally-2kn}
Let $n\geq 2$. Let $\G$ be a connected locally $2\K_n$ graph. Let $\Sigma=C(\G)$ and $\a\in V(\Sigma)$. Then $\G$ is $3$-CSH but not $3$-CH if and only if one of the following holds:
\begin{enumerate}
  \item [{\rm (1)}]\  $\Sigma$ is a tetravalent $3$-arc-regular graph and $\Aut(\Sigma)_\a\cong \A_4\rtimes C_3$;
  \item [{\rm (2)}]\  $\Sigma$ is a pentavalent $3$-arc-regular graph  and $\Aut(\Sigma)_\a\cong \F(20)\times C_4$;
  \item [{\rm (3)}]\  $\Sigma$ is a $3$-arc-transitive graph of valency $8$ and $\Aut(\Sigma)_\a^{\Sigma(\a)}\cong C_2^3\rtimes (C_{7}\rtimes C_{t})$ with $t=1$ or $3$;
  \item [{\rm (4)}]\  $\Sigma$ is a $3$-arc-transitive graph of valency $32$ and $\Aut(\Sigma)_\a^{\Sigma(\a)}\cong C_2^5\rtimes (C_{31}\rtimes C_{5})$;
  \item [{\rm (5)}]\ $\Sigma$ is a $3$-arc-transitive graph of valency $q+1$ and $\Aut(\Sigma)_\a^{\Sigma(\a)}\cong\PSL(2,q).\lg\eta\rg$, where $q$ is an odd prime power such that $q\equiv-1\ (\mod 4)$ and $\eta$ is a field automorphism of $\GF(q)$;
  \item [{\rm (6)}]\  $\Sigma$ is a pentavalent symmetric graph of type ${\mathcal{Q}_2^6}$ and $\Aut(\Sigma)_\a\cong \F(20)\times C_2$;

  \item [{\rm (7)}]\  $\Sigma$ is a trivalent symmetric graph of type $2^2$ and $\Aut(\Sigma)_\a\cong \S_3$.
  \end{enumerate}
\end{theorem}

\f\demo From Lemmas~\ref{lem:triangle-tran}--\ref{lem:not2gt} we immediately obtain the sufficiency. So we only need to prove the necessity.

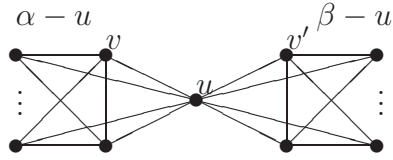
\begin{figure}[ht]
\begin{center}
\unitlength 4mm
\begin{picture}(12,6)

\put(0, 5.5){$\a-{u}$} \put(10, 5.5){$\b-{u}$}

\put(6, 3){\circle*{0.4}} \put(6, 3.2){$u$}

\put(6,3){\line(2, 1){3}} \put(6,3){\line(2, -1){3}}
\put(6,3){\line(-2, 1){3}} \put(6,3){\line(-2, -1){3}}

\put(6,3){\line(4, 1){6}} \put(6,3){\line(4, -1){6}}
\put(6,3){\line(-4, 1){6}} \put(6,3){\line(-4, -1){6}}

\put(3, 4.5){\circle*{0.4}} \put(3, 4.7){$v$}
\put(3, 1.5){\circle*{0.4}} 
\put(0, 4.5){\circle*{0.4}} 
\put(0, 1.5){\circle*{0.4}} 
\put(0, 2.5){$\vdots$} \put(12, 2.5){$\vdots$}

\put(0,1.5){\line(1, 1){3}}
\put(0,1.5){\line(1, 0){3}}
\put(3,4.5){\line(0, -1){3}} \put(3,1.5){\line(-1, 1){3}}
\put(3,4.5){\line(-1, 0){3}}

\put(9, 4.5){\circle*{0.4}} \put(9, 4.7){$v'$}
\put(9, 1.5){\circle*{0.4}} 
\put(12, 4.5){\circle*{0.4}} 
\put(12, 1.5){\circle*{0.4}} 

\put(9,1.5){\line(0, 1){3}}
\put(9,1.5){\line(1, 1){3}}
\put(9,1.5){\line(1, 0){3}}
\put(12,1.5){\line(-1, 1){3}}
\put(12,4.5){\line(-1, 0){3}}
\end{picture}
\end{center}\vspace{-.5cm}
\caption{The subgraph of $\G$ induced by $\{u\}\cup\G(u)$} \label{fig-4}
\end{figure}

Suppose that $\G$ is a locally $2\K_n$ graph which is $3$-CSH. By Lemma~\ref{lem:3CSH-2arc}, $\Sigma$ is $2$-arc-transitive. Take $u\in V(\G)$.  In what follows, we always assume that $\a$ and $\b$ are the two maximal cliques of $\G$ containing $u$ (see Figure~(\ref{fig-4})). Let $A=\Aut(\G)$. Then every element of $A_u$ either fixes or interchanges $\a-\{u\}$ and $\b-\{u\}$ as $\G(u)=(\a-\{u\})\cup(\b-\{u\})$. Clearly, $\{\a, \b\}$ is an edge of $\Sigma$ and $\a\cap\b=\{u\}$. Furthermore, for any $v\in\a-\{u\}$, we have $\a-\{u,v\}=\G(u)\cap\G(v)$. By Lemma~\ref{3-csh-stabilizer}, $A_u$ is both vertex-transitive and edge-transitive on $[\G(u)]^c$.
We shall divide the proof into the following two cases:

\medskip
\f{\bf Case~1.}\ $A_u$ is arc-transitive on $[\G(u)]^c$.

By Lemma~\ref{3-csh-stabilizer}~(5), $\G$ is $2$-geodesic transitive, and then by Proposition~\ref{prop:jinwei-line}, $\Sigma$ is $3$-arc-transitive. If $A_u$ is arc-transitive on $[\G(u)]$, then by \cite[Proposition~2.1]{Li-Zhou}, $\G$ is $3$-CH, a contradiction.


Thus, $A_u$ is not arc-transitive on $[\G(u)]$. By Lemma~\ref{3-csh-stabilizer}~(3), $A_{vu}$ has $s$ orbits of equal size on $\a-\{u,v\}$ with $s=1, 2, 3$ or $6$, where $v\in\a-\{u\}$. Note that $\a-\{u,v\}=\G(u)\cap\G(v)$. If $s=1$, then $A_{vu}$ is transitive on $\G(u)\cap\G(v)$, and then $A_u$ would be arc-transitive on $[\G(u)]$ as $A_u$ is transitive on $\G(u)$, a contradiction. Thus, $s\neq 1$.
Since $A_u$ is transitive on $\G(u)=(\a-\{u\})\cup(\b-\{u\})$, it follows that $A_{\a-\{u\}}$ is transitive on $\a-\{u\}$. This implies that the permutation group $(A_\a^\a)_u$ is either regular or $\frac{3}{2}$-transitive on $\a-\{u\}$.

Assume first that $(A_\a^\a)_u$ is regular on $\a-\{u\}$. Then $A_{uv}$ fixes every vertex in $\a-\{u\}$. Since $A_{vu}$ has $s$ orbits of equal size on $\a-\{u,v\}$ with $s=2, 3$ or $6$, it follows that $|\a-\{u,v\}|=2, 3$ or $6$, and hence $|\a-\{u\}|=3, 4$ or $7$. If $|\a-\{u\}|=3$, then we have $A_\a^\a\cong\A_4$ and $|\a|=4$. So $\Sigma$ has valency $4$. Since $\Sigma$ is $3$-arc-transitive,  by \cite[Theorem~4]{Potocnik}, we have that $A_\a\cong \A_4\times C_3$ and so $\Sigma$ is a tetravalent $3$-arc-regular graph. This implies part (1). If $|\a-\{u\}|=4$, then we have $|\a|=5$ and $A_\a^\a\cong\F(20)$. So $\Sigma$ has valency $5$. Since $\Sigma$ is $3$-arc-transitive, by \cite[Theorem~4.1]{Zhou-Feng-2010}, we have $A_\a\cong \F(20)\times C_4$, and hence $\Sigma$ is $3$-arc-regular. Then part (2) happens. If $|\a-\{u\}|=7$, then we have $|\a|=8$ and $A_\a^\a\cong C_2^3\rtimes C_7$. So $\Sigma$ has valency $8$ and part (3) happens.

Now assume that $(A_{\a}^\a)_u$ is $\frac{3}{2}$-transitive on $\a-\{u\}$. Then $(A_{\a}^\a)_u$ has rank $3, 4$ or $7$ on $\a-\{u\}$ as $A_{uv}$ has $s$ orbits of equal size on $\a-\{u,v\}$ with $s=2, 3$ or $6$. Clearly, $A_\a^\a$ is a $2$-transitive permutation group on $\a$. By Lemma~\ref{lem:2-tran+impri}, $(A_\a^\a)_u$ is primitive on $\a-\{u\}$ and Lemma~\ref{lem:2-tran+impri}~(1) or (2) happens. If Lemma~\ref{lem:2-tran+impri}~(2) happens, then either part (3) or part (4) holds. If Lemma~\ref{lem:2-tran+impri}~(1) happens, by Lemmas~\ref{lem:triangle-tran-1}--\ref{lem:triangle-tran}, we obtain part (5).



\medskip
\f{\bf Case~2.} $A_u$ is not arc-transitive on $[\G(u)]^c$.

In this case, $\G$ is not $2$-geodesic transitive, and then by Proposition~\ref{prop:jinwei-line}, $\Sigma$ is $2$-arc-transitive but not $3$-arc-transitive. Recall that $A_u$ is both vertex-transitive and edge-transitive on $[\G(u)]^c$. It follows that $A_u$ is half-arc-transitive on $[\G(u)]^c$.  Since $[\G(u)]\cong 2\K_n$, one has $[\G(u)]^c\cong\K_{n,n}$, and so $n$ is even. Let $B=\a-\{u\}$ and $C=\b-\{u\}$. Then $|B|=|C|=n$. Since $A_u$ is transitive on $\G(u)$, $(A_u)_B$ is transitive on $B$. Since $n-1=|B|-1$ is odd, by Lemma~\ref{3-csh-stabilizer}~(3), $A_{uv}$ has $s$ orbits of equal size on $B-\{v\}=\G(u)\cap \G(v)$ with $s=1$ or $3$, where $v\in B$. In particular, the permutation group $(A_\a^\a)_u$ is either regular or $\frac{3}{2}$-transitive on $\a-\{u\}$.

Suppose first that $A_{uv}$ has $3$ orbits of equal size on $B-\{v\}$. Clearly, $A_\a^\a$ is a $2$-transitive permutation group on $\a$. By Lemma~\ref{lem:2-tran+impri}, $(A_\a^\a)_u$ can not be a $\frac{3}{2}$-transitive permutation group on $B=\a-\{u\}$ of rank $4$. So $(A_\a^\a)_u$ is regular on $B$, and hence all orbits of $A_{uv}$ on $B-\{v\}$ have size $1$. Then $|(A_\a^\a)_u|=|B|=4$ and $|\a|=5$. It then follows that $A_\a^\a\cong\F(20)$ and $\Sigma$ is a pentavalent $2$-arc-transitive graph. Since $A_{u}$ is half-arc-transitive on $[\G(u)]^c\cong\K_{4,4}$, one has $A_u^{\G(u)}$ is a $2$-group. Note that $A_u$ is just the stabilizer of the edge $\{\a,\b\}$ of $\Sigma$ in $A$. By \cite[Theorem~1.2]{Morgan}, we see that $A_\a\cong \F(20)\times C_2$ and $A_{\{\a,\b\}}\cong {\rm M}_{16}$. We obtain part~(6).

Now suppose that $A_{uv}$ has only one orbit on $B-\{v\}$. Then $(A_\a^\a)_u$ is $2$-transitive on $B=\a-\{u\}$. Since $A_u$ is transitive on $\G(u)$, it follows that $(A_u^{\G(u)})_B$ is $2$-transitive on $B$ and $(A_u^{\G(u)})_C$ is also $2$-transitive on $C$. If $(A_u^{\G(u)})_B$ is not faithful on $B$, then  the kernel of $(A_u^{\G(u)})_B$ on $B$ would be transitive on $C$. However, this is impossible since $A_u$ is half-arc-transitive on $[\G(u)]^c$. Thus, $(A_u^{\G(u)})_B$ is faithful on $B$. Similarly,  $(A_u^{\G(u)})_C$ is faithful on $C$. Note that $(A_u^{\G(u)})_B=(A_u^{\G(u)})_C$. By \cite{Cameron-2tr}, $(A_u^{\G(u)})_B$ is either affine or almost simple.




We now claim that the actions of $(A_u^{\G(u)})_B$ on $B$ and $C$ are equivalent. Take $v\in B$. Since $A_u$ is half-arc-transitive on $[\G(u)]^c\cong\K_{n,n}$, $A_{uv}$ has two orbits on $C$ of equal size. In particular, $|B|=|C|=n$ is even.

Assume first that $(A_u^{\G(u)})_B$ is affine. Then $\soc((A_u^{\G(u)})_B)$ is elementary abelian of order $|B|$. Since $|B|$ is even, we have $\soc((A_u^{\G(u)})_B)\cong C_2^r$ for some integer $r>0$. If the actions of $(A_u^{\G(u)})_B$ on $B$ and $C$ are not equivalent, then by inspecting the affine $2$-transitive permutation groups (see \cite[Table~7.3]{Cameron-Book}), we conclude that $r>2$ and one of the following may happens:
\begin{enumerate}
  \item [{\rm (i)}]\  $A_{uv}^{\G(u)} \leq \G{\rm L}(1,2^r)\cong C_{2^r-1}\rtimes C_r$;
  \item [{\rm (ii)}]\ $\SL(d,q)\unlhd A_{uv}^{\G(u)} \leq \Gamma {\rm L}(d,q)(2^r=q^d, d\geq 2)$;
  \item [{\rm (iii)}]\ ${\rm Sp}(d,q)\unlhd A_{uv}^{\G(u)}\leq (\mz_{q-1}\circ{\rm Sp}(2d,q)).(\mz_{(2,q-1)}\times\mz_{r/2d}) (q^{2d}=2^r, d\geq 2)$;
  \item [{\rm (iv)}]\ ${\rm G}_2(q)\unlhd A_{uv}^{\G(u)}\leq(\mz_{q-1}\circ{\rm G}_2(q)).\mz_{r/6} (q^{6}=2^r)$;
  \item [{\rm (v)}]\ $A_{uv}^{\G(u)}=\A_6$ and $r=4$;
  \item [{\rm (vi)}]\ $A_{uv}^{\G(u)}=\PSU(3,3)$ and $r=6$.
\end{enumerate}
In case (i), as $A_{uv}$ has two orbits on $C$ of length $2^{r-1}$, it follows that $2^{r-1}\mid |A_{uv}^{\G(u)}|$, and so $2^{r-1}\mid r(2^r-1)$. It implies that $2^{r-1}\mid r$. This, however, is impossible because $r>2$.
In cases (ii), (iii) or (iv), the center $Z$ of $A_{uv}^{\G(u)}$ has order dividing $q-1$. Take $w\in C$. As $A_{uv}$ has two orbits on $C$ of length $2^{r-1}$, it follows that $|w^Z|$ is divisor of gcd$(q-1, 2^{r-1})=1$, and hence $|w^Z|=1$. This implies that $Z$ fixes every vertex in $C$. Since $(A_u^{\G(u)})_B$ is faithful on $C$, one has $Z=1$. Since $r>2$, $A_{uv}^{\G(u)}$ is almost simple. Now we conclude that $\soc(A_{uv}^{\G(u)})$ is one of the following groups:
\begin{equation}\label{simplegroups}
\PSL(d,q) (2^r=q^d), \soc({\rm PSp}(d, q)) (q^{2d}=2^r), {\rm G}_2(q) (q^6=2^r), \A_6, \PSU(3,3).
\end{equation}
On the other hand, since $(A_u^{\G(u)})_B$ is faithful on $C$ and $A_{uv}$ has two orbits on $C$ of length $2^{r-1}$, it follows that $\soc(A_{uv}^{\G(u)})$ has an orbit on $C$ of length dividing $2^{r-1}$. This implies that $\soc(A_{uv}^{\G(u)})$ has a maximal subgroup of index $2^t$ with $1<t\leq r-1$. By \cite{Guralnick} (or \cite[Theorem~2.2]{LMP2009}), we have $\soc(A_{uv}^{\G(u)})\cong\A_{2^t}$ or $\PSL(2,\ell)$ with $\ell$ a prime and $\ell+1=2^t$. Since $\soc(A_{uv}^{\G(u)})$ is also one of the groups in (\ref{simplegroups}), by the Classification Theorem for Finite Simple Groups (see for example \cite[p.3]{Wilson}) we conclude that either $r=4$ and $\SL(4, 2)=\soc(A_{uv}^{\G(u)})\cong\A_{8}$, or $r=3$ and $\SL(3, 2)=\soc(A_{uv}^{\G(u)})\cong\PSL(2,7)$. For the former, we have $r=4$ and $A_{u}^{\G(u)}=\AGL(4,2)$. By {\tt Magma}~\cite{BCP}, all subgroups of $\AGL(4,2)$ isomorphic to $\A_8$ are conjugate, and so the actions of $(A_u^{\G(u)})_B$ on $B$ and $C$ are equivalent, a contradiction. For the latter, we have $r=3$ and $\soc(A_{uv}^{\G(u)})\cong\PSL(2,7)$, but $\PSL(2,7)$ does not have a subgroup of index no more than $2^2$, a contradiction.

Assume now that $(A_u^{\G(u)})_B$ is almost simple. If the actions of $(A_u^{\G(u)})_B$ on $B$ and $C$ are not equivalent, then by checking \cite[Theorem~5.3]{Cameron-2tr},  one of the following holds:
\begin{enumerate}
\item [{\rm (a)}] $\soc((A_u^{\G(u)})_B)=\A_6\,$ and $|B|=6$;
\item [{\rm (b)}] $\soc((A_u^{\G(u)})_B)=\PSL(d,q) (d>2)\,$ and $|B|=\frac{q^d-1}{q-1}$;
\item [{\rm (c)}] $(A_u^{\G(u)})_B={\rm M}_{12}$ (Mathieu)\, and $|B|=12$;
\item [{\rm (d)}] $(A_u^{\G(u)})_B={\rm HS} $(Higman-Sims)\, and $|B|=176$.
\end{enumerate}

In cases (a) and (c), by {\tt Magma}~\cite{BCP}, we can obtain that $(A_u^{\G(u)})_v$ is transitive on $C$, a contradiction.
In case (d), by {\tt Magma}~\cite{BCP}, $(A_u^{\G(u)})_v$ has exactly two orbits on $C$ with size $50$ and $126$, respectively. This is contrary to the fact that the two orbits of $(A_u^{\G(u)})_v$ on $C$ has equal size.

In case (b), we have $\soc((A_u^{\G(u)})_B)=\PSL(d,q) (d>2)$ and $|B|=\frac{q^d-1}{q-1}$. Here we may assume that $B$ and $C$ are the set of points and the set of hyperplanes of the projective space
$PG(d\!-\!1, q)$, respectively. Then the hyperplanes containing $v$ form an orbit $C_1$ of $\soc((A_u^{\G(u)})_B)_v$ on $C$, while the hyperplanes not containing $v$ form another orbit $C_2$ of $\soc((A_u^{\G(u)})_B)_v$ on $C$. Then $|C_1|=\frac{q^{d-1}-1}{q-1}$, and $|C_2|=\frac{q^d-q^{d-1}}{q-1}$. Since the two orbits of $\soc((A_u^{\G(u)})_B)_v$ on $C$ have equal size, we have $|C_1|=|C_2|$, and hence $q^{d-1}-1=q^d-q^{d-1}$, namely, $q^{d}+1=2q^{d-1}$. This, however, is impossible.

By now, we have shown that the actions of $(A_u^{\G(u)})_B$ on $B$ and $C$ are equivalent. So $(A_u^{\G(u)})_v$ also fixes at least one vertex in $C$. Again, since $(A_u^{\G(u)})_v$ has exactly two orbits of equal size on $C$, we have $|C|=2$. Then $\G$ has valency $4$ and $\Sigma$ has valency $3$. Since $\Sigma$ is $2$- but not $3$-arc-transitive, one has $|A_u|=4$. Again, since $A_u$ acts half-arc-transitively on $[\G(u)]^c\cong\K_{2,2}$, we must have $A_u\cong C_4$, and so $\Sigma$ is a trivalent symmetric graph of type $2^2$. Part~(7) happens.\hfill\qed

\section{Proof of Theorem~\ref{th:sol-3csh-graphs}}
The goal of this section is to characterize solvable $3$-CSH but not $3$-CH graphs and prove Theorem~\ref{th:sol-3csh-graphs}.
We first give several lemmas about arc-transitive graphs.
Let $\G$ be a $(G,s)$-arc-transitive graph with $G\leq\Aut(\G)$ and $s\geq 2$, and let $N$ be a normal subgroup of $G$. The {\em
quotient graph} $\G_N$ of $\G$ relative to $N$ is defined as the graph
with vertices the orbits of $N$ on $V(\G)$ and with two different orbits
adjacent if there exists an edge in $\G$ between the vertices lying in
those two orbits. If $\G_N$ and $\G$ have the same valency, then we say that $\G$ is a {\em normal cover} of $\G_N$. 
In view of \cite[Theorem~4.1]{Praeger1993} or \cite[Lemma~2.5]{Li-Pan}, we have the following.

\begin{lem}\label{lem:quot}
Let $\G$ be a connected $(G,2)$-arc-transitive graph with $G\leq\Aut(\G)$. Suppose that $N\unlhd G$ has at least three orbits on $V(\G)$. Then
\begin{enumerate}
  \item [{\rm (1)}]\ $N$ acts semiregularly on $V(\G)$ and $\G$ is a normal cover of $\G_N$.
  \item [{\rm (2)}]\ $N$ is the kernel of $G$ acting on $V(\G_N)$, $G/N\leq\Aut(\G_N)$ and $\G_N$ is $(G/N, 2)$-arc-transitive.
\end{enumerate}
\end{lem}

\begin{lem}\label{lem:abelian-edge-regular}
Let $\G$ be a connected graph of valency $k>2$. Suppose that $G\leq\Aut(\G)$ is abelian and acts regularly on the edge set of $\G$. Then $\G\cong\K_{k,k}$.
\end{lem}

\f\demo Take an edge $\{u,v\}\in E(\G)$. Since $G$ is regular on the edge set of $\G$, one has $E(\G)=\{\{u^g,v^g\}\mid g\in G\}$ and
$|G|=|E(\G)|$. Then $V(\G)=u^G\cup v^G$. Since $k>2$, one has $|V(\G)|<|E(\G)|$, and so $u^G, v^G$ are two distinct orbits of $G$ on $V(\G)$. Furthermore, $\G$ is a bipartite graph with two bi-parts $u^G$ and $v^G$. By the edge-transitivity of $G$ on $\G$, we have $G_u$ is transitive on $\G(u)$. Clearly, $v\in \G(u)$. For any $w\in \G(v)-\{u\}$, since $G$ is abelian, we have $G_u=G_w$, and so $w$ is adjacent to all vertices in $\G(u)$ by the transitivity of $G_u$ on $\G(u)$. By the arbitrariness of $w$, we see that the subgraph induced by $\G(u)\cup \G(v)$ is isomorphic to $\K_{k,k}$. Since $\G$ is connected, one has $\G\cong\K_{k,k}$.\hfill\qed

\begin{lem}\label{lem:line-cay}
Let $\G$ be a connected graph of valency $k$. Suppose that $G\leq\Aut(\G)$ is regular on $E(\G)$ and intransitive on $V(\G)$. Let $\{u,v\}$ be an edge of $\G$. Then the following hold.
\begin{enumerate}
  \item [{\rm (1)}]\ $G=\lg G_u, G_v\rg$.
  \item [{\rm (2)}]\ $G_u\cap G_v=1$.
  \item [{\rm (3)}]\ $|G_u|=|G_v|=k$, and $G_u$ and $G_v$ act regularly on $\G(u)$ and $\G(v)$, respectively.
  \item [{\rm (4)}]\ The line graph of $\G$ is isomorphic to the Cayley graph $\Cay(G, S)$ with $S=(G_u\cup G_v)-\{1\}$.
\end{enumerate}
\end{lem}

\f\demo Since $G$ is regular on $E(\G)$ but intransitive on $V(\G)$, the connectedness of $\G$ implies that $G=\lg G_u, G_v\rg$. This proves part (1).

Again since $G$ is regular on $E(\G)$, one has $G_{\{u,v\}}=1$ and $|G|=|E(\G)|$. Moreover, $V(G)=u^G\cup v^G$. Since $G$ is intransitive on $V(\G)$, one has $u^G\cap v^G=\emptyset$. It follows that $G_u\cap G_v=G_{\{u, v\}}=1$, proving part (2).

Since $u^G\cap v^G=\emptyset$, neither $u^G$ nor $v^G$ contains an edge of $\G$. It follows that $\G(u)\subseteq v^G$ and $\G(v)\subseteq u^G$. As $G$ is transitive on $E(G)$, $G_u$ and $G_v$ are transitive on $\G(u)$ and $\G(v)$, respectively. Moreover, $|E(G)|=|u^G||\G(u)|=|v^G||\G(v)|$. Note that $|G|=|u^G||G_u|=|v^G||G_v|$. Since $|G|=|E(\G)|$, it follows that $|G_u|=|G_v|=|\G(u)|=|\G(v)|=k$. Then $G_u$ and $G_v$ act regularly on $\G(u)$ and $\G(v)$, respectively. We obtain part (3).

As $G$ is regular on $E(\G)$, the line graph, say $\Sigma$, of $\G$ is a Cayley graph on $G$. Note that $V(\Sigma)=E(\G)$ and the set of vertices of $\Sigma$ adjacent to $\{u, v\}$ is
\[F=\{\{u, x\}, \{v, y\}\mid v\neq x\in\G(u), u\neq y\in\G(v)\}.\]
Clearly, $|F|=|\G(u)-\{v\}|+|\G(v)-\{u\}|=2k-2$.
Let $S=\{g\in G\mid \{u,v\}^g\in F\}$. Then $\Sigma\cong\Cay(G, S)$ and $|S|=|F|$. For any $g\in G_u-\{1\}$, we have $\{u,v\}^g=\{u, v^g\}$ with $v\neq v^g\in\G(u)$, and hence $\{u,v\}^g\in F$. It follows that $g\in S$ and so $G_u-\{1\}\subseteq S$. Similarly, we have $G_v-\{1\}\subseteq S$. Since $|(G_u\cup G_v)-1|=|F|$, one has $S=(G_u\cup G_v)-1$. This proves part (4).
\hfill\qed

\begin{lem}\label{th:sol-s-arc-tran}
Let $\G$ be a connected $(G,s)$-arc-transitive graph, where $s\geq 2$ and $G\leq\Aut(\G)$ is solvable. Then either
\begin{enumerate}
  \item [{\rm (1)}]\ $G$ has a normal subgroup which is semiregular on $V(\G)$ with at most $2$ orbits and for each vertex $v$ of $\G$, $G_v$ acts faithfully on $\G(v)$; or
  \item [{\rm (2)}]\ $\G$ is a normal cover of $\K_{p^n,p^n}$ with $p$ a prime and $n$ a positive integer, and $G$ has a normal subgroup, say $M$ such that the following hold:
\begin{enumerate}
  \item [{\rm (i)}]\ $M$ is regular on $E(\G)$ and intransitive on $V(\G)$, $M=\lg M_u, M_w\rg$, $M_u\cap M_w=1$ and $M_u\cong M_w\cong C_p^n$, where $\{u,w\}\in E(\G)$; and
  \item [{\rm (ii)}]\ the line graph of $\G$ is isomorphic to the Cayley graph $\Cay(M, S)$ with $S=(M_u\cup M_w)-\{1\}$.
\end{enumerate}

\end{enumerate}
\end{lem}

\f\demo Let $N\unlhd G$ be maximal subject to the condition that $N$ has at
least three orbits on $V(\G)$. Let $\G_N$ be the quotient graph of $\G$ relative to $N$. Since $G$ is $2$-arc-transitive on $\G$, by Lemma~\ref{lem:quot}, $\G$ is a normal cover of $\G_N$, and $N$ is semiregular on $V(\G)$ and $N$ is the kernel of $G$ acting on $V(\G_N)$. Furthermore, $G/N$ is a group of automorphisms of $\G_N$ acting transitively on the $2$-arcs of $\G_N$. Let $M/N$ be a minimal normal subgroup of $G/N$. The solvability of $G$ implies that $M/N\cong C_p^r$ with $p$ a prime and $r$ a positive integer. By the maximality of $N$, either $M/N$ is transitive on $V(\G_N)$ or $M/N$ has two orbits on $V(\G_N)$. If the former happens, then $M/N$ is regular on $V(\G_N)$, and then by the semiregularity of $N$ on $V(\G)$, $M$ is a normal subgroup of $G$ acting regularly on $V(\G)$. So $\G$ is a Cayley graph on $M$. Take an arbitrary $u\in V(\G)$, and let $S=\{g\in M\mid \{u, u^g\}\in E(\G)\}$. Then $\G\cong\Cay(M, S)$. Without loss of generality, we may let $\G=\Cay(M, S)$. Then $\G(1)=S$. Since $M\unlhd G$, one has $G_1\leq\Aut(M, S)=\{\a\in\Aut(M)\mid S^\a=S\}$ (see \cite{Godsil1981}). Since $\G$ is connected, one has $M=\lg S\rg$. This implies that the vertex stabilizer $G_1$ acts faithfully on $\G(1)=S$. Since $\G$ is vertex-transitive, for each vertex $v$ of $\G$, $G_v$ acts faithfully on $\G(v)$, as claimed in part (1).

Now let $M/N$ have two orbits on $V(\G_N)$. Then $M$ has two orbits, say $U$ and $W$, on $V(\G)$. So $U,W$ are blocks of imprimitivity of $G$ on $V(\G)$. Since $G$ is $2$-arc-transitive on $\G$, $U$ and $W$ contain no edges of $\G$, and so $\G$ is bipartite with $U$ and $W$ as its two bi-parts. If $M/N$ is semiregular on $V(\G_N)$, then since $N$ is semiregular on $V(\G)$, $M$ is also semiregular on $V(\G)$. By \cite[Lemma~2.4]{Li-2005-PAMS}, we see that for each vertex $v$ of $\G$, $G_v$ acts faithfully on $\G(v)$, as claimed in part (1).


Now suppose that $M/N$ is not semiregular on $V(\G_N)$. Let $\{u,w\}$ be an edge of $\G$ such that $u\in U$ and $w\in W$. Let $B=u^N$ and $C=w^N$. Then $M_B=M_uN$ and $M_C=M_wN$. Recall that $N$ is semiregular on $V(\G)$. It follows that $M_u\cong M_B/N$ and $M_w\cong M_C/N$. As $M/N$ is not semiregular on $V(\G_N)$, $M_B/N$ is a non-trivial normal subgroup of $G_B/N$. Since $G/N$ is $2$-arc-transitive on $\G_N$, $M_B/N$ is transitive on the neighbors of $B$ in $\G_N$. It follows that $M/N$ is transitive on the edges of $\G_N$. Then $M/N$ is regular on the edges of $\G_N$ since $M/N\cong C_p^r$. By Lemma~\ref{lem:abelian-edge-regular}, we have $\G_N\cong\K_{p^n, p^n}$, and $M/N=M_B/N\times M_C/N\cong C_p^{n}\times C_p^n$ with $r=2n$. Furthermore, $M_u\cong M_B/N\cong C_p^n$ and $M_w\cong M_C/N\cong C_p^n$. Since $\G$ is a normal cover of $\G_N$ and since $M/N$ is regular on $E(\G_N)$, it implies that $M$ is regular on $E(\G)$. Since $M/N$ has two orbits on $V(\G_N)$, it follows that $M$ is intransitive on $V(\G)$. Applying Lemma~\ref{lem:line-cay} to $\G$ and $M$ we can obtain (i) and (ii) of part (2).
\hfill\qed

Now we are ready to prove Theorem~\ref{th:sol-3csh-graphs}. We first prove the necessity of Theorem~\ref{th:sol-3csh-graphs} in the following lemma.

\begin{lem}\label{lem:sol-3csh-graphs}
Let $n\geq 2$ and let $\G$ be a solvable locally $2\K_n$ graph. If $\G$ is $3$-CSH, then
$\G$ is isomorphic to an arc-transitive normal Cayley graph $\Cay(H,S)$ on a group $H$ such that the following hold:
\begin{enumerate}
  \item [{\rm (a)}]\ $H$ has two subgroups $A,B$ such that $H=\lg A,B\rg$, $A\cong B\cong C_p^f$, $A\cap B=1$ and $S=(A\cup B)-\{1\}$; and
  \item [{\rm (b)}]\ if $\G$ is not $3$-CH, then one of $(1)-(6)$ of Theorem~{\rm \ref{th:sol-3csh-graphs}}~$(b)$ holds.
\end{enumerate}
\end{lem}

\f\demo  Suppose that $\G$ is $3$-CSH. If $\G$ is $3$-CH, then by Theorem~\ref{mainth}, $C(\G)$ is $3$-arc-transitive and locally $3$-transitive, and if $\G$ is not $3$-CH, then by Theorem~\ref{mainth:locally-2kn}, $C(\G)$ is a $2$-arc-transitive graph satisfying the conditions in one of parts (1)--(4), (6) and (7) of Theorem~\ref{mainth:locally-2kn}. Let $\Sigma=C(\G)$. By Lemma~\ref{lem:clique-line}, $\G$ is isomorphic to the line graph of $\Sigma$. For convenience, in the following, we shall identify $\G$ with the line graph of $\Sigma$. Due to Lemma~\ref{lem:clique-line}, we may also view $\Aut(\G)$ as the full automorphism group of $\Sigma$. Take a vertex $u$ of $\Sigma$ and take $v\in \Sigma(u)$.

We begin by proving that $\Aut(\G)_u$ is not faithful on $\Sigma(u)$ if $\Sigma$ is not the case in part (7) of Theorem~\ref{mainth:locally-2kn}. First, if $\Sigma$ is $3$-arc-transitive and locally $3$-transitive, then since $\Aut(\G)$ is solvable, it follows that $\Aut(\G)_u^{\Sigma(u)}$ is a solvable $3$-transitive permutation group on $\Sigma(u)$. By checking the list of finite affine 2-transitive permutation groups obtained by Hering (see for example \cite[Table 7.3]{Cameron-Book}), we see that either $\Sigma$ has valency $3$ and $\Aut(\G)_u^{\Sigma(u)}\cong\S_3$, or $\Sigma$ has valency $4$ and $\Aut(\G)_u^{\Sigma(u)}\cong\S_4$. Since $\Sigma$ is $3$-arc-transitive, by Theorem~\ref{th:3-arc-tran-stabi}, $\Aut(\G)_u$ is not faithful on $\Sigma(u)$, as claimed.

If $\Sigma$ is a graph in part (2) of Theorem~\ref{mainth:locally-2kn}, then $\Sigma$ is a pentavalent $3$-arc-regular graph, and then by \cite[Table~2]{Morgan}, we see that $\Aut(\G)_u\cong \F(20)\times C_4$, $\Aut(\G)_{\{u,v\}}\cong C_4\wr C_2$ and $\Aut(\G)_u^{[1]}\cong C_4$. Similarly, if $\Sigma$ is a graph in part (1) of Theorem~\ref{mainth:locally-2kn}, then $\Sigma$ is a tetravalent $3$-arc-regular graph, and by \cite[Theorem~4]{Potocnik}, we see that $\Aut(\G)_u\cong \A_{4}\times C_3$, $\Aut(\G)_{\{u,v\}}\cong C_4\wr C_2$ and $\Aut(\G)_u^{[1]}\cong C_3$. Moreover, if $\Sigma$ is a graph in part (6) of Theorem~\ref{mainth:locally-2kn}, then it is a pentavalent $3$-arc-transitive graph of type $\mathcal{Q}_2^6$. By \cite[Table~1]{Morgan}, we obtain that $\Aut(\G)_u\cong \F(20)\times C_2$, $\Aut(\G)_{\{u,v\}}\cong {\rm M}_{16}$ and $\Aut(\G)_u^{[1]}\cong C_2$.
For the graph $\Sigma$ in part (3) of Theorem~\ref{mainth:locally-2kn}, since $\Sigma$ is $3$-arc-transitive, one has $8\cdot 7^2\mid |\Aut(\G)_u|$, and in view of \cite[Theorem~2.1]{Li-S-Song-2015-JA}, we see that
\[\begin{array}{c}
(C_2^3\rtimes C_7)\times C_7\leq \Aut(\G)_u\leq (C_2^3\rtimes (C_7\rtimes C_3))\times (C_7\rtimes C_3),\\
(C_7\times C_7)\rtimes C_2\leq \Aut(\G)_{\{u,v\}}\leq ((C_7\rtimes C_3)\times (C_7\rtimes C_3))\rtimes C_2.
\end{array}\]
Similarly, for the graph $\Sigma$ in part (4) of Theorem~\ref{mainth:locally-2kn}, we have $2^5\cdot 31^2\mid |\Aut(\G)_u|$, and by \cite[Theorem~2.1]{Li-S-Song-2015-JA}, we have
\[\begin{array}{c}
(C_2^5\rtimes (C_{31}\rtimes C_5))\times C_{31}\leq \Aut(\G)_u\leq (C_2^5\rtimes (C_{31}\rtimes C_5))\times (C_{31}\rtimes C_5),\\
(C_{31}\times C_{31})\rtimes C_{10}\leq \Aut(\G)_{\{u,v\}}\leq ((C_{31}\rtimes C_5)\times (C_{31}\rtimes C_{5}))\rtimes C_2.
\end{array}\]

So far we have shown that if $\Sigma$ is not the case in part (7) of Theorem~\ref{mainth:locally-2kn}, then $\Sigma$ is a $2$-arc-transitive graph and $\Aut(\G)_u$ is not faithful on $\Sigma(u)$. Since $\Aut(\G)$ is solvable, applying Lemma~\ref{th:sol-s-arc-tran} to $\Aut(\G)$ we see that Lemma~\ref{th:sol-s-arc-tran}~(2) happens. It follows that $\G$ is a normal Cayley graph $\Cay(H,S)$ on a group $H$, where $S=(A\cup B)-\{1\}$ with $A=H_u$ and $B=H_v$. Furthermore, $H=\lg A,B\rg$, $A\cap B=1$, and $A\cong B\cong C_p^f$, where $p$ is a prime and $p^f$ is just the valency of $\Sigma$. So $H$ and $S$ satisfies the condition in part (a).

If part~(7) of Theorem~\ref{mainth:locally-2kn} happens, then $\Sigma$ is a trivalent symmetric graph of type $2^2$. Then $\Aut(\G)_u\cong \S_3$, and for any $v\in\Sigma(u)$, the edge stabilizer $\Aut(\G)_{\{u,v\}}$ is isomorphic to $C_4$. By \cite[Corollary~1.2]{FengLiZhou}, $\Aut(\G)$ has a normal subgroup, say $N$, such that the quotient graph $\Sigma_N$ of $\Sigma$ relative to $N$ is isomorphic to $\K_{3,3}$. So $\Aut(\G)/N$ is $2$-arc-transitive on $\Sigma_N\cong\K_{3,3}$. Let $P/N$ be the Sylow $3$-subgroup of $\Aut(\G)/N$. Then $P/N\cong C_3^2$, which is normal in $\Aut(\G)/N$ and is regular on the edges of $\Sigma_N\cong\K_{3,3}$. Then $P$ is regular on $E(\Sigma)$ but intransitive on $V(\Sigma)$. Furthermore, $P_u\cong P_v\cong C_3$ and $P=\lg P_u, P_v\rg$. Now applying Lemma~\ref{lem:line-cay} to $P$ and observing that $P\unlhd\Aut(\G)$, we see that $\G$ is a normal Cayley graph $\Cay(H,S)$ on $H$, where $H=P$, $A=P_u$, $B=P_v$, $H=\lg A, B\rg$ and $S=(A\cup B)-\{1\}$. So $H$ and $S$ satisfy the condition in part (a).

As a conclusion, $\G$ is a normal Cayley graph $\Cay(H, S)$ on a group $H$, where $S=(A\cup B)-\{1\}$, $A$ and $B$ are subgroups of $H$ such that $H=\lg A, B\rg$, $A\cap B=1$ and $A\cong B\cong\mz_p^f$. Note that $\Aut(H,S)$ is just the stabilizer of an edge of $\Sigma$ in $\Aut(\G)$. From the argument in the above paragraphs, one may see that if $\G$ is not $3$-CH, then one of $(1)-(6)$ of Theorem~{\rm \ref{th:sol-3csh-graphs}}~$(b)$ holds. \hfill\qed

Finally, we prove the sufficiency of Theorem~\ref{th:sol-3csh-graphs}. The main thing that we need to prove is the Cayley graph $\Cay(H, S)$ satisfying the conditions in part (a) and part (b) of Theorem~\ref{th:sol-3csh-graphs} is 3-CSH but not $3$-CH. Actually, we can give more information about the symmetry of $\Cay(H, S)$. This is done in the following lemma.

\begin{lem}\label{lem:testfornormal2geod-t}
Let $p$ be a prime and $f$ be a positive integer. Let $H$ be a group having two subgroups $A,B$ such that $H=\lg A,B\rg$, $A\cong B\cong C_p^f$, and $A\cap B=1$. Let $\G=\Cay(H,S)$ with $S=(A\cup B)-\{1\}$, and let $\Sigma=C(\G)$. Then $\Sigma$ has valency $p^f$, and $R(H)$ is regular on $E(\Sigma)$ and intransitive on $V(\Sigma)$. Moreover, if $\Aut(\G)$ is solvable and $\G$ satisfies the conditions in one of $(1)$--$(6)$ of Theorem~\ref{th:sol-3csh-graphs}~$(b)$, then the following hold.
    \begin{enumerate}
   \item [{\rm (1)}]\ If $(p,f)=(2,2)$ and $\Aut(H,S)\cong C_3\wr C_2$, then $\Sigma$ is a tetravalent $3$-arc-regular graph, and $\G$ is a $6$-valent $3$-CSH and $2$-geodesic-transitive but not $3$-CH graph;
   \item [{\rm (2)}]\ If $(p,f)=(5,1)$ and $\Aut(H,S)\cong C_4\wr C_2$, then $\Sigma$ is a pentavalent $3$-arc-regular graph, and $\G$ is an $8$-valent $3$-CSH and $2$-geodesic-transitive but not $3$-CH graph;
   \item [{\rm (3)}]\ If $(p,f)=(2,3)$ and $C_7\wr C_2\leq \Aut(H,S)\leq (C_7\rtimes C_3)\wr C_2$, then $\Sigma$ is a $3$-arc-transitive graph of valency $8$, and $\G$ is a $14$-valent $3$-CSH and $2$-geodesic-transitive but not $3$-CH graph;
   \item [{\rm (4)}]\ If $(p,f)=(2,5)$ and $(C_{31}\times C_{31})\rtimes C_{10}\leq \Aut(H,S)\leq (C_{31}\rtimes C_5)\wr C_2$, then $\Sigma$ is a $32$-valent $3$-arc-transitive graph, and $\G$ is a $62$-valent $3$-CSH and $2$-geodesic-transitive but not $3$-CH graph;
  \item [{\rm (5)}]\ If $(p,f)=(5,1)$ and $\Aut(H,S)\cong {\rm M}_{16}$, then $\Sigma$ is a pentavalent symmetric graph of type ${\mathcal{Q}_2^6}$, and $\G$ is an $8$-valent $3$-CSH but not $2$-geodesic-transitive graph;
  \item [{\rm (6)}]\ If $(p,f)=(3,1)$ and $\Aut(H,S)\cong C_4$, then $\Sigma$ is a trivalent symmetric graph of type $2^2$, and $\G$ is a tetravalent $3$-CSH but not $2$-geodesic-transitive graph.
  \end{enumerate}
\end{lem}

\f\demo Since $\G=\Cay(H,S)$ with $S=(A\cup B)-\{1\}$, it follows that both $A$ and $B$ induce two subgraphs of $\G$ isomorphic to $\K_{p^f}$ due to $A\cong B\cong C_p^f$. As $S=(A\cup B)-\{1\}$ and $A\cap B=1$, $[A]$ and $[B]$ are two maximal cliques of $\G$ containing $1$. Since $R(H)$ acts transitively on $V(\G)=H$ by right multiplication, $\{[Ah], [Bh]\mid h\in H\}$ is the set of all maximal cliques of $\G$. For arbitrary $h\in H$, we have $Ah\cap Bh=\{h\}$ and $\G(h)=(Ah\cup Bh)-\{h\}$ as $A\cap B=\{1\}$ and $S=(A\cup B)-\{1\}$. This implies that every vertex $h$ of $\G$ is contained in exactly two maximal cliques, namely, $[Ah]$ and $[Bh]$. Let $\Sigma=C(\G)$ be the clique graph of $\G$. Then $V(\Sigma)=\{Ah, Bh\mid h\in H\}$ and $E(\Sigma)=\{\{Ah, Bh\}\mid h\in H\}$. It follows that $\Sigma$ is a bipartite graph with two bi-parts $B_0=\{Ah\mid h\in H\}$ and $B_1=\{Bh\mid h\in H\}$. Clearly, $R(H)$ induces an action on $V(\Sigma)$ with two orbits $B_0$ and $B_1$. Since $A\cap B=1$, one has $R(H)_A\cap R(H)_B=\{R(g)\mid g\in A\cap B\}=1$. This implies that $R(H)$ acts faithfully on $V(\Sigma)$. Since $\Sigma(A)=\{Ba\mid a\in A\}$, it follows that $\Sigma$ has valency $|A|=p^f$ and $R(H)_A$ acts regularly on $\Sigma(A)$ as $R(H)_A\cap R(H)_B=1$. Thus, $R(H)$ is regular on $E(\Sigma)$ but not transitive on $V(\Sigma)$. By Lemma~\ref{lem:clique-line}, $\G$ is isomorphic to the line graph $L(\Sigma)$ of $\Sigma$. This proves the first half of our lemma.


In what follows, we shall prove the second half of this lemma. For convenience, we identify $\G$ with the line graph of $\Sigma$. Then $\Aut(H, S)$ fixes an edge, say $\{u,v\}$, of $\Sigma$. Let $N=R(H)\rtimes\Aut(H,S)$. Since $H=\lg A, B\rg$ and $S=(A\cup B)-\{1\}$, $\Aut(H,S)$ acts faithfully on $S$. Note that both $[A]$ and $[B]$ are two maximal cliques of $\G$. Since $A\cap B=1$ and $\Aut(H, S)$ fixes $1$, it implies that every element of $\Aut(H,S)$ either setwise fixes or interchanges the subsets $A-\{1\}$ and $B-\{1\}$ of $S$. Then $\Aut(H, S)$ acts on $\{A, B\}$, and let $K$ be the kernel of this action. Then $\Aut(H,S)\cong\Aut(H,S)^S\leq(K^{A}\times K^{B})\rtimes C_2$. Since $\Aut(H,S)$ acts faithfully on $(A\cup B)-\{1\}$, it implies that $K^A$ is faithful on $A$. Since $K\leq\Aut(H,S)\leq\Aut(H)$ and $A\leq H$, one has $K^A\leq\Aut(A)$. Similarly, $K^B\leq\Aut(B)$. It follows that $\Aut(H,S)\cong\Aut(H,S)^S\leq (\Aut(A)\times\Aut(B))\rtimes C_2$. 

We first deal with the case when $(p,f)=(2,f)$ with $f=2,3$ or $5$. Note that if $f=2$, then $\Aut(H,S)\cong C_3\wr C_2$, and if $f=3$, then $C_7\wr C_2\leq \Aut(H,S)$, and if $f=5$, then $(C_{31}\times C_{31})\rtimes C_{10}\leq \Aut(H,S)\leq (C_{31}\rtimes C_5)\wr C_2$. It follows that in these three cases, $N_u$ is not faithful on $\Sigma(u)$ and that $N$ is $3$-arc-transitive on $\Sigma$. Furthermore, if $f=2$, then $N_u^{\Sigma(u)}\cong \A_4$, if $f=3$, then $N_u^{\Sigma(u)}\cong C_2^3\rtimes (C_7\rtimes C_s)$ with $s=1$ or $3$, and if $f=5$, then $N_u^{\Sigma(u)}\cong C_2^5\rtimes (C_{31}\rtimes C_5)$. To prove that part (1), (3) and (4) hold, by Lemma~\ref{lem:triangle-tran}, it suffices to prove that $N=\Aut(\Sigma)$, namely, $R(H)\unlhd \Aut(\Sigma)$. Since $\Aut(\G)$ is solvable, by Lemma~\ref{th:sol-s-arc-tran}, the following statements are true:
\begin{enumerate}
  \item [{\rm (a)}]\ $\Aut(\Sigma)$ has a normal subgroup, say $T$, which is regular on $E(\Sigma)$ and intransitive on $V(\Sigma)$,
  \item [{\rm (b)}]\ the edge-stabilizer $\Aut(\Sigma)_{\{u,v\}}\leq (C_{2^f-1}\rtimes C_f)\wr C_2$, and,
  \item [{\rm (c)}]\ $\Sigma$ is a bipartite graph with two bi-parts of size $\frac{|T|}{2^f}=\frac{|R(H)|}{2^f}$.
\end{enumerate}
To complete the proof, it suffices to prove that $R(H)=T$. Suppose by way of contradiction that $R(H)\neq T$. Noticing that both $R(H)$ and $T$ are edge- but not vertex-transitive on $\Sigma$, $R(H)T$ is edge- but not vertex-transitive on $\Sigma$. Let $D$ be the stabilizer of the edge $\{u, v\}$ of $\Sigma$ in $R(H)T$. Then $R(H)T=T\rtimes D$, and $D$ also fixes the arc $(u,v)$. Note that arc-stabilizer $\Aut(\Sigma)_{(u,v)}\leq (C_{2^f-1}\rtimes C_f)\times (C_{2^f-1}\rtimes C_f)$. So $1\neq D\lesssim K\times K$, where $K=\S_3$ if $f=2$, or $K=C_7\rtimes C_3$ if $f=3$ or $K=C_{31}\rtimes C_5$ if $f=5$. Let $M=R(H)\cap T$. Then $M\unlhd N$. Suppose that $M$ is semiregular on $V(\Sigma)$. Then the number of orbits of $M$ on $V(\Sigma)$ is equal to $2\cdot\frac{|R(H)|}{2^f\cdot |M|}=\frac{1}{2^{f-1}}\cdot\frac{|R(H)T|}{|T|}=\frac{|D|}{2^{f-1}}$. This means $4\mid |D|$, forcing $f=2$ and $D\leq \S_3\times\S_3$. It follows that the number of orbits of $M$ on $V(\Sigma)$ is equal to $2$, $6$ or $18$. If $M$ has two orbits on $V(\Sigma)$, then by \cite[Lemma~2]{Li-2005-PAMS}, $N_u$ is faithful on $\Sigma(u)$, which is impossible since $N$ is $3$-arc-transitive on $\Sigma$ and $\Sigma$ has valency $4$. If $M$ has $6$ or $18$ orbits on $V(\Sigma)$, then by Lemma~\ref{lem:quot}, the quotient graph of $\Sigma$ relative to $M$ is a tetravalent $2$-arc-transitive graph of order $6$ or $18$. However, by \cite[Table~3]{Potocnik}, no such graph exists. Thus, $M$ is not semiregular on $V(\Sigma)$, and so $M_u\neq 1$. Since $N$ is $3$-arc-transitive on $\Sigma$, we see that $N_u$ is $2$-transitive on $\Sigma(u)$, and since $1\neq M_u\unlhd N_u$, it implies that $M_u$ is transitive on $\Sigma(u)$. So $M$ is edge-transitive on $\Sigma$. It follows that $M=R(H)=T$ since $R(H)$ and $T$ are regular on $E(\Sigma)$. This completes the proof of part (1), (3) and (4).

Suppose that $(p,f)=(5,1)$. Then $\Sigma$ has valency $5$, and $|A-\{1\}|=|B-\{1\}|=4$. Recall that $\Aut(H,S)\cong\Aut(H,S)^S\leq(\Aut(A)\times\Aut(B))\rtimes C_2$. It follows that if $\Aut(H,S)\cong C_4\wr C_2$ or ${\rm M}_{16}$, then $\Aut(H,S)$ acts transitively on $S$, and then $\G$ is arc-transitive. By Lemma~\ref{prop:jinwei-line}, $\Sigma$ is $2$-arc-transitive. Assume $\Aut(H,S)\cong C_4\wr C_2$. By \cite[Theorem~1.2]{Morgan}, $N_u\cong \F(20)\times C_4$ and $N_u^{[1]}\cong C_4$.  Since $\Aut(\Sigma)$ is solvable, by \cite[Theorem~4.1]{Zhou-Feng-2010} or \cite[Theorem~1.2]{Morgan}, we have $\Aut(\Sigma)=N$, and so $\Sigma$ is $3$-arc-regular. Again, by Lemma~\ref{lem:3-ar-2gt}, $\G$ is $2$-geodesic-transitive but not $3$-CH. This implies part (2). Assume now $\Aut(H,S)\cong {\rm M}_{16}$. Then by \cite[Theorem~1.2]{Morgan}, we have $N_u\cong \F(20)\times C_2$ and $N_u^{[1]}\cong C_2$. Since $\Aut(\G)$ is solvable, by \cite[Theorem~4.1]{Zhou-Feng-2010} or \cite[Theorem~1.2]{Morgan},  one has $\Aut(\Sigma)_u\leq \F(20)\times C_4$, and so $\Aut(\Sigma)_u$ has a unique Sylow $5$-subgroup. Moreover, by Theorem~\ref{th:sol-s-arc-tran}, $\Aut(\Sigma)$ contains a normal subgroup, say $T$, such that $T$ is regular on $E(\Sigma)$. It follows that $T_u$ is a Sylow $5$-subgroup of $\Aut(\Sigma)_u$. Clearly, $R(H)_u$ is also a Sylow $5$-subgroup of $\Aut(\Sigma)_u$, so $R(H)_u=T_u$. Similarly, $R(H)_v=T_v$. It follows that $R(H)=\lg R(H)_u,R(H)_v\rg=\lg T_u,T_v\rg=T$, and hence $R(H)\unlhd \Aut(\Sigma)$. Thus, $N=R(H)\rtimes\Aut(H,S)=\Aut(\G)$, and hence $\Sigma$ is a pentavalent symmetric graph of type $\mathcal{Q}_2^6$. By Lemma~\ref{lem:not2gt}, $\G$ is $3$-CSH but not $2$-geodesic-transitive. This proves part (5).

Finally, suppose that $(p,f)=(3,1)$ and $\Aut(H,S)\cong C_{4}$. Then $\Sigma$ is a trivalent graph. Since $|S|=4$, $\Aut(H,S)$ acts regularly on $S$, and so $N$ is arc-transitive on $\G$. It follows that $N$ is $2$-arc-transitive on $\Sigma$. Since $\Aut(H,S)$ is the stabilizer of the edge $\{u,v\}$ in $N$, one has $N_u\cong \S_{3}$. By \cite[Theorem~5.1]{ConderNedela}, $\Sigma$ is $2$- or $3$-arc-regular. So, $\Aut(\Sigma)_u\leq \S_3\times C_2$, and so $\Aut(\Sigma)_u$ has a unique Sylow $3$-subgroup. By Theorem~\ref{th:sol-s-arc-tran}, the solvability of $\Aut(\G)$ implies that $\Aut(\Sigma)$ has a normal subgroup, say $T$, acting regularly on $E(\Sigma)$. It follows that $T_u$ is a Sylow $3$-subgroup of $\Aut(\Sigma)_u$. Clearly, $R(H)_u$ is also a Sylow $3$-subgroup of $\Aut(\Sigma)_u$, so $R(H)_u=T_u$. Similarly, $R(H)_v=T_v$. It follows that $R(H)=\lg R(H)_u,R(H)_v\rg=\lg T_u,T_v\rg=T$, and hence $R(H)\unlhd \Aut(\Sigma)$. Thus, $N=R(H)\rtimes\Aut(H,S)=\Aut(\G)$, and so $\Sigma$ is trivalent symmetric graph of type $2^2$. By Lemma~\ref{lem:not2gt}, $\G$ is $3$-CSH but not $2$-geodesic-transitive. This proves part (6). \hfill\qed



\section{Examples of $3$-CSH but not $3$-CH graphs}\label{subsec1}

We begin by constructing a graph satisfying the condition (1) of Theorem~\ref{th:sol-3csh-graphs}~(b).\medskip

\f{\bf Construction~I}\
Let $\H=\lg a,b,c,d,e,f,g,h\rg$ be a group with the following relations:
\[\begin{array}{l}
a^2=b^2=c^4=d^4=e^2=f^2=g^4=h^4=1,[a,b]=[e,f]=1,\\
c=[a,e], d=[a,f], g=[b,e], h=[b,f],[c,d]^2=[c,g]^2=[c,h]^2=[d,g]^2=[d,h]^2=[g,h]^2=1, \\
\end{array}\]
\[\begin{array}{l}
[a,[c,d]]=[b,[c,d]]=[e,[c,d]]=[f,[c,d]]=[a,[c,g]]=[b,[c,g]]=[e,[c,g]]=[f,[c,g]]=1,\\

[a,[c,h]]=[b,[c,h]]=[e,[c,h]]=[f,[c,h]]=[a,[g,d]]=[b,[g,d]]=[e,[g,d]]=[f,[g,d]]=1,\\

[a,[h,d]]=[b,[h,d]]=[e,[h,d]]=[f,[h,d]]=[a,[g,h]]=[b,[g,h]]=[e,[g,h]]=[f,[g,h]]=1,\\
\end{array}\]
\[\begin{array}{l}
[[a,e],[b,f]][f,[a,[b,e]]]=[[b,e],[ab,f]][f,[b,[ab,e]]]=[[ab,e],[a,f]][f,[ab,[a,e]]]=1,\\

[[a,f],[b,ef]][ef,[a,[b,f]]]=[[b,f],[ab,ef]][ef,[b,[ab,f]]]=[[ab,f],[a,ef]][ef,[ab,[a,f]]]=1,\\

[[a,ef],[b,e]][e,[a,[b,ef]]]=[[b,ef],[ab,e]][e,[b,[ab,ef]]]=[[ab,ef],[a,e]][e,[ab,[a,ef]]]=1,\\
\end{array}\]
\[\begin{array}{l}
[[e,a],[f,b]][b,[e,[f,a]]]=[[f,a],[ef,b]][b,[f,[ef,a]]]=[[ef,a],[e,b]][b,[ef,[e,a]]]=1,\\

[[e,b],[f,ab]][ab,[e,[f,b]]]=[[f,b],[ef,ab]][ab,[f,[ef,b]]]=[[ef,b],[e,ab]][ab,[ef,[e,b]]]=1,\\

[[e,ab],[f,a]][a,[e,[f,ab]]]=[[f,ab],[ef,a]][a,[f,[ef,ab]]]=[[ef,ab],[e,a]][a,[ef,[e,ab]]]=1.
\end{array}
\]
Let $\Delta=\Cay(\H,S)$ with $S=\{a,b, ab, e,f, ef\}$.

\begin{lem}
The group $\H$ has order $2^{17}$, and the graph $\Delta$ is a normal Cayley graph on $\H$ with $\Aut(\H,S)\cong C_3^2\rtimes C_2$. In particular, $\Delta$ is $2$-geodesic-transitive and $3$-CSH but not $3$-CH, and $C(\Delta)$ is a tetravalent $3$-arc-regular Cayley graph of order $2^{16}$.
\end{lem}

\f\demo We shall make use of {\tt Magma~\cite{BCP}} in the proof and see Appendix~1 for the programs. Let $A=\lg a, b\rg $ and $B=\lg e,f\rg$. By using the {\tt pQuotient} command in {\tt Magma}~\cite{BCP}, we obtain that $|\H|=2^{17}$, $A\cong B\cong C_2^2$, $H=\lg A, B\rg$ and $A\cap B=1$. Then $S=(A\cup B)-1$ and the subgraphs of $\Delta$ induced by $A$ and $B$, respectively, are two maximal cliques $\K_{4}$. This implies that $\Aut(\H,S)$ acts on $\{A, B\}$. Moreover, $\Aut(\H,S)$ acts faithfully on $A\cup B$ since $\H=\lg A,B\rg$. It follows that $\Aut(\H, S)\cong\Aut(\H,S)^S\leq (\Aut(A)\times\Aut(B))\rtimes C_2\cong\S_3\wr C_2$.

By using the {\tt hom} command in {\tt Magma}~\cite{BCP}, we see that both the map $a\mapsto b, b\mapsto a, e\mapsto e, f\mapsto f$ and the map $a\mapsto b, b\mapsto a, e\mapsto f, f\mapsto e$ do not induce automorphisms of $\H$, but each of the following maps induces an automorphism of $\H$:
\[
\begin{array}{l}
\a: a\mapsto b, b\mapsto ab, e\mapsto e, f\mapsto f, \\

\b: a\mapsto a, b\mapsto b, e\mapsto f, f\mapsto ef,\\


\g: a\mapsto e, b\mapsto f, e\mapsto a, f\mapsto b.

\end{array}
\]
It follows that $\Aut(\H,S)=(\lg\a\rg\times\lg\b\rg)\rtimes\lg\g\rg\cong C_3^2\rtimes C_2$.

Let $\Sigma=C(\Delta)$. Recall that $R(H)$ acts on $V(\Delta)=H$ by right multiplication. Since $A\cap B=1$, for all $h\in H$, we have $Ah\cap Bh=h$ and $\Delta(h)=(Ah\cup Bh)-h$ as $\Delta(1)=S=(A\cup B)-1$. This implies that $\Sigma$ has vertex set $V(\Sigma)=\{Ah, Bh\mid h\in \H\}$, and edge set $\{\{Ah, Bh\}\mid h\in H\}$. Let $V_0=\{Ah\mid h\in H\}$ and $V_1=\{Bh\mid h\in H\}$. Then $|V_0|=|V_1|=2^{15}$, and so $|V(\Sigma)|=2^{16}$. Moreover, $V_0$ and $V_1$ are two orbits of $R(H)$ on $V(\Sigma)$. Note that the neighborhood of $A$ in $\Sigma$ is $\Sigma(A)=\{Ba\mid a\in A\}$ while the neighborhood of $B$ in $\Sigma$ is $\Sigma(B)=\{Ab\mid b\in B\}$. As $A\cong B\cong C_2^2$, it follows that $\Sigma$ has valency $4$. Let $K=\lg ae, bf, \H'\rg$. By {\tt Magma}~\cite{BCP}, we have $A\cap K=B\cap K=1$. It follows that $|H: K|=2^4$ and hence $|H|=|V_0|=|V_1|$. Since $R(H)_A=\{R(a)\mid a\in A\}$ and $R(H)_B=\{R(b)\mid b\in B\}$, it implies that $R(K)=\{R(g)\mid g\in K\}$ acts semiregularly on $V(\Sigma)$ with two orbits $V_0$ and $V_1$. Notice that $\g$ swaps $a$ and $e$, and swaps $b$ and $f$. It follows that $\g$ centralizes $R(ae)$ and $R(bf)$. As $\H'\unlhd \H$, we have $\g$ normalizes $R(K)$, implying that $|\lg R(K),\g\rg|=|V(\Sigma)|$. Note that $\g$ also swaps $A$ and $B$. This implies that $\lg R(K),\g\rg$ is transitive and so regular on $V(\Sigma)$. It follows that $\Sigma$ is a Cayley graph.

Since $\Sigma$ is a tetravalent graph, the stabilizer of any vertex of $\Sigma$ in $\Aut(\Sigma)$ is a $\{2,3\}$-group, and hence $\Aut(\Sigma)$ is also a $\{2,3\}$-group since $\Sigma$ has order $2^{16}$. It follows that $\Aut(\Sigma)$ is solvable. By Lemma~\ref{lem:testfornormal2geod-t}~(1), $\Delta$ is $2$-geodesic-transitive and $3$-CSH, but not $3$-CH, and $\Sigma$ is a tetravalent $3$-arc-regular Cayley graph.\hfill\qed

The following proposition proves that there exist infinitely many solvable tetravalent $3$-arc-regular graphs.

\begin{prop}\label{family-1}
There exist infinitely many solvable tetravalent $3$-arc-regular graphs.
\end{prop}

\f\demo By \cite[Theorem~2.11~(1)]{Du-Kwak-Xu-laa-2003}, for all primes $p>2^2\cdot 3^4$ there exists a connected $(X, 3)$-arc-regular graph $\Pi$ with $X\leq\Aut(\Pi)$ satisfying the following conditions:
\begin{enumerate}
  \item [{\rm (1)}]\  $X$ has a normal subgroup $N\cong\C_p^{\b(\G)}$, where $\b(\G)=|E(C(\Delta))|-|V(C(\Delta))|+1$ is the Betti number of $C(\Delta)$ and $\Delta$ is the graph in Construction~I;
  \item [{\rm (2)}]\  the norma quotient $\Pi_N\cong C(\Delta)$ and $\Pi$ is a normal cover of $\Pi_N$;
  \item [{\rm (3)}]\  $X/N\cong\Aut(C(\Delta))$.
\end{enumerate}
We claim that $\Pi$ is $3$-arc-regular. Suppose on the contrary that $\Pi$ is not $3$-arc-regular. Then $\Aut(\Pi)>X$. Since $\Pi$ is a tetravalent $2$-arc-transitive graph, by \cite[Theorem~4]{Potocnik} we have $|\Aut(\Pi)|\mid 2^4\cdot 3^6\cdot |V(C(\Delta))|$. So $N$ is a Sylow $p$-subgroup of $\Aut(\Pi)$. Since $X$ is $3$-arc-regular on $\Pi$, one has $|X|=2^2\cdot 3^2\cdot|V(C(\Delta))|$. Consequently, $|\Aut(\Pi): X|\mid 2^2\cdot 3^4$. Since $N\unlhd X$, one has $|\Aut(\Pi): N_{\Aut(\Pi)}(N)|\mid 2^2\cdot 3^4$, and since $p>2^2\cdot 3^4$, by Sylow's theorem, we have $\Aut(\Pi)=N_{\Aut(\Pi)}(N)$, and so $N\unlhd \Aut(\Pi)$. By Lemma~\ref{lem:quot}, we would have $\Aut(\Pi)/N\leq\Aut(C(\Delta))$, which is impossible because $C(\Delta)$ is $3$-arc-regular.
\hfill\qed

\medskip
\f{\bf Remark on Proposition~\ref{family-1}.}\ From \cite[Theorem~1.1~\&~Corollary~1.2]{LiLingMa} one may deduce that every tetravalent $3$-arc-regular Cayley graph is a normal cover of a Cayley graph on one of the following groups: $C_3^{11}\rtimes (C_2^{12}.{\rm M}_{11})$, $\S_{35}$ and $\A_{35}$. This, however, is not true by Proposition~\ref{family-1}.\medskip

Next we give two graphs satisfying the condition (2) of Theorem~\ref{th:sol-3csh-graphs}~(b), of which the first one appeared in \cite[Remark~4.2]{Zhou-EJC-2021}.

\begin{exam}\label{exam:5-val-1}
{\rm Let $\M=\lg a,b,c\ |\ a^5=b^5=c^5=1, c=[a,b], [a,c]=[b,c]=1\rg$. Let $\Theta=\Cay(\M,S)$ with $S=\{a,a^2,a^3,a^4, b,b^2,b^3,b^4\}$.
By {\tt Magma}~{\rm \cite{BCP}}, $\M$ has order $5^{3}$, and the graph $\Theta$ is a normal Cayley graph on $\M$ with $\Aut(\M,S)\cong C_4^2\rtimes C_2$. By Lemma~{\rm \ref{lem:testfornormal2geod-t}~(2)}, $\Theta$ is $2$-geodesic-transitive but not $3$-CH, and  $C(\Theta)$ is a pentavalent $3$-arc-regular graph of order $2\cdot 5^{2}$. (See Appendix~3 for the {\tt Magma} programs used in this example.)}
\end{exam}

\begin{exam}\label{exam:5-val-2}
{\rm Let $\T=\lg a,b,c,d,e\ |\ a^5=b^5=c^5=d^5=e^5=1, c=[a,b], d=[a,c], e=[b,c], [a,d]=[b,d]=[a,e]=[b,e]=1\rg$. Let $\Phi=\Cay(\T,S)$ with $S=\{a,a^2,a^3,a^4, b,b^2,b^3,b^4\}$. By {\tt Magma}~{\rm \cite{BCP}}, $\T$ has order $5^{5}$, and the graph $\Phi$ is a normal Cayley graph on $\T$ with $\Aut(\T,S)\cong C_4^2\rtimes C_2$. By Lemma~{\rm \ref{lem:testfornormal2geod-t}~(2)}, $\Phi$ is $2$-geodesic-transitive but not $3$-CH, and $C(\Phi)$ is a pentavalent $3$-arc-regular graph of order $2\cdot 5^{4}$. (See Appendix~4 for the {\tt Magma} programs used in this example.)}
\end{exam}

The following example gives a graph satisfying the condition (3) of Theorem~\ref{th:sol-3csh-graphs}~(b).

\begin{exam}\label{exam:cond-3}
{\rm Let $\mathcal{G}=\lg a,b,c,e,f,g,x,y,z\rg$ be a group with the following relations:
\[\begin{array}{l}
a^2=b^2=c^2=e^2=f^2=g^2=x^2=y^2=z^2=1, [g,a]=x, [g,b]=y, [g,c]=z,\\

[a,b]=[a,c]=[b,c]=[e,f]=[e,g]=[f,g]=1,\\

[e,a]=xyz, [e,b]=xz, [e,c]=x, [f,a]=xz, [f,b]=x, [f,c]=y.
\end{array}\]
Let $\Theta=\Cay(\mathcal{G},S)$, where $S=A\cup B-\{1\}$, $A=\lg a,b,c\rg$ and $B=\lg e,f,g\rg$.
By {\tt Magma}~{\rm \cite{BCP}}, $\mathcal{G}$ has order $2^{9}$, and the graph $\Theta$ is a normal Cayley graph on $\mathcal{G}$ with $\Aut(\mathcal{G},S)\cong (C_7\times C_7)\rtimes C_6$. By Lemma~{\rm \ref{lem:testfornormal2geod-t}~(3)},  $\Theta$ is $2$-geodesic-transitive and $3$-CSH but not $3$-CH, and $C(\Theta)$ is a $3$-arc-transitive graph. (See Appendix~5 for the {\tt Magma} programs used in this example.)}
\end{exam}

The following example gives a graph satisfying the condition (4) of Theorem~\ref{th:sol-3csh-graphs}~(b).

\begin{exam}\label{exam:cond-4}
{\rm Let $\L=\lg a,b,c,d,e,u,v,x,y,z,f,g,h,i,j\rg$ be a group with the following relations:
\[\begin{array}{l}
a^2=b^2=c^2=d^2=e^2=f^2=g^2=h^2=i^2=j^2=u^2=v^2=x^2=y^2=z^2=1, \\

[a,b]=[a,c]=[a,d]=[a,e]=[b,c]=[b,d]=[b,e]=[c,d]=[c,e]=[d,e]=1,\\

[u,v]=[u,x]=[u,y]=[u,z]=[v,x]=[v,y]=[v,z]=[x,y]=[x,z]=[y,z]=1,\\

[a,z]=f,[b,z]=g,[c,z]=h,[d,z]=i,[e,z]=j,\\

[a,u]=fgh,[b,u]=ghi,[c,u]=hij,[d,u]=fghj,[e,u]=f,\\

[a,v]=ghi,[b,v]=hij,[c,v]=fghj,[d,v]=f,[e,v]=g\\

[a,x]=hij,[b,x]=fghj,[c,x]=f,[d,x]=g,[e,x]=h,\\

[a,y]=fghj,[b,y]=f,[c,y]=g,[d,y]=h,[e,y]=i.

\end{array}\]
Let $\Pi=\Cay(\L,S)$, where $S=A\cup B-\{1\}$, $A=\lg a,b,c,d,e\rg$ and $B=\lg u,v,x,y,z\rg$.
By {\tt Magma}~{\rm \cite{BCP}}, $\L$ has order $2^{15}$, and the graph $\Pi$ is a normal Cayley graph on $\L$ with $\Aut(\L,S)\cong (C_{31}\times C_{31})\rtimes C_{10}$. By Lemma~{\rm \ref{lem:testfornormal2geod-t}~(4)},  $\Pi$ is $2$-geodesic-transitive and $3$-CSH but not $3$-CH, and $C(\Pi)$ is a $3$-arc-transitive graph. (See Appendix~6 for the {\tt Magma} programs used in this example.)}
\end{exam}

Now we give a graph satisfying  the condition (5) of Theorem~\ref{th:sol-3csh-graphs}~(b).

\begin{exam}\label{exam:val-5-m16}
{\rm Let $\N=\lg a,b,c,d,e,f,g,h,k\rg$ be a group with the following relations:
\[\begin{array}{l}
a^5=b^5=c^5=d^5=e^5=f^5=g^5=h^5=k^5=1,\\

c=[a,b],d=[a,c],e=[b,c],[d,e]=1, [a,d]=f, [b,d]=g, [a,e]=h, [b,e]=k,\\

[a,f]=[a,g]=[a,h]=[a,k]=[b,f]=[b,g]=[b,h]=[b,k]=1, f=k^{-2}, g=h^{-2}.
\end{array}\]
Let $\Lambda=\Cay(\N,S)$ with $S=\{a,a^2,a^3,a^4, b,b^2,b^3,b^4\}$.
By {\tt Magma}~{\rm \cite{BCP}}, $\N$ has order $5^{6}$, and the graph $\Lambda$ is a normal Cayley graph on $\N$ with $\Aut(\N,S)\cong {\rm M}_{16}$. By Lemma~{\rm \ref{lem:testfornormal2geod-t}~(5)},  $\Lambda$ is $3$-CSH but not $2$-geodesic-transitive, and $C(\Lambda)$ is a pentavalent symmetric graph of type ${\mathcal{Q}_2^6}$. (See Appendix~7 for the {\tt Magma} programs used in this example.)}
\end{exam}

\f{\bf Remark on Examples~\ref{exam:5-val-1}--\ref{exam:val-5-m16}.}\  {\rm With a similar argument as in the proof of Proposition~\ref{family-1}, by using the graphs given in Examples~\ref{exam:5-val-1}--\ref{exam:val-5-m16}, one can see that there exist infinitely many graphs satisfying the conditions (2)--(5) of Theorem~\ref{th:sol-3csh-graphs}~(b).}\medskip

Finally, we construct a family of graphs satisfying the condition (6) of Theorem~\ref{th:sol-3csh-graphs}~(b).\medskip

\f{\bf Construction~II}\
Let $n\geq 2$ be an integer, and let $\R=\lg a,b\rg$ be a finite $3$-group with the following relations:
\[
\begin{array}{l}
a^3=b^3=c^{3^n}=d^{3^n}=e^{3^n}=f^{3^n}=g^3=h^3=1, c=[a,b], d=[b,a^2], e=[a^2,b^2], f=[b^2,a], \\

[c,d]=c^{-3^{n-1}}d^{3^{n-1}}, [c,f]=c^{-3^{n-1}}f^{3^{n-1}},[d,e]=d^{-3^{n-1}}e^{3^{n-1}},[e,f]=e^{-3^{n-1}}f^{3^{n-1}},\\

[d^{3^{n-1}},c]=[d^{3^{n-1}},e]=[d^{3^{n-1}},f]=1,[e^{3^{n-1}},c]=[e^{3^{n-1}},d]=[e^{3^{n-1}},f]=1,\\

f^{3^{n-1}}=c^{3^{n-1}}d^{-3^{n-1}}e^{3^{n-1}}, [c^{3^{n-1}},d]=[c^{3^{n-1}},e]=[c^{3^{n-1}},f]=1,\\

g=[c,e], h=[d,f], [g,a]=[g,b]=[h,a]=[h,b]=1, \\

h=c^{3^{n-1}}d^{3^{n-1}}e^{3^{n-1}}, g^{-1}=d^{3^{n-1}}e^{3^{n-1}}f^{3^{n-1}}. \\
\end{array}
\]
Let $\Upsilon=\Cay(\R,S)$ with $S=\{a,a^2,b,b^2\}$.\medskip

\begin{lem}\label{lem:cubic2-2type}
The group $\R$ has order $3^{4n+1}$, and the graph $\Upsilon$ is a normal Cayley graph on $\R$ with $\Aut(\R,S)\cong C_4$. In particular, $\Upsilon$ is $3$-CSH but not $2$-geodesic-transitive, and $C(\Upsilon)$ is a trivalent symmetric graph of type $2^2$.
\end{lem}

\f\demo Let $D=\lg c,d,e,f\rg$. By a direct calculation, we obtain the following relations:
\[c^a=d^{-1}e^{-1}, d^a=c, e^a=f, f^a=f^{-1}e^{-1}, c^b=c^{-1}f^{-1},d^b=e^{-1}d^{-1},e^b=d, f^b=c.\]
Since $\R$ is generated by $a,b$, one has $D\unlhd\R$ and $\R/D=\lg aD, bD\rg$. Since $c=[a,b]\in D$, one has $[aD,bD]=[a,b]D=D$ and hence $\R/D$ is abelian. As both $a$ and $b$ have order $3$, $\R/D=\lg aD\rg\times\lg bD\rg\cong C_3\times C_3$. It follows that $D=\Phi(\R)=\R'$.

We now show the following two claims.

\medskip
\f{\bf Claim~1}\ Let $a_1=b$ and $b_1=a^2$. Then $a_1$ and $b_1$ have the same relations as do $a$ and $b$. \medskip

Let $c_1=[a_1,b_1]$, $d_1=[b_1,a_1^2]$, $e_1=[a_1^2,b_1^2]$ and $f_1=[b_1^2,a_1]$. Then $c_1=d,$ $d_1=e,$ $e_1=f$ and $f_1=c$.
So $c_1^{3^n}=d_1^{3^n}=e_1^{3^n}=f_1^{3^n}=1$. From the following relations
\[
\begin{array}{l}
f^{3^{n-1}}=c^{3^{n-1}}d^{-3^{n-1}}e^{3^{n-1}},\\

[c^{3^{n-1}},d]=[c^{3^{n-1}},e]=[c^{3^{n-1}},f]=1,\\

[d^{3^{n-1}},c]=[d^{3^{n-1}},e]=[d^{3^{n-1}},f]=1,\\

[e^{3^{n-1}},c]=[e^{3^{n-1}},d]=[e^{3^{n-1}},f]=1,
\end{array}
\]
we know that $c^{3^{n-1}}, d^{3^{n-1}}, e^{3^{n-1}},f^{3^{n-1}}$ are in the center of $D$.
So
\[f_1^{3^{n-1}}=c^{3^{n-1}}=f^{3^{n-1}}d^{3^{n-1}}e^{-3^{n-1}}=c_1^{3^{n-1}}d_1^{-3^{n-1}}e_1^{3^{n-1}},\]
and $c_1^{3^{n-1}}, d_1^{3^{n-1}}, e_1^{3^{n-1}}$ are in the center of $D$.

Let $g_1=[c_1,e_1]$ and $h_1=[d_1,f_1]$. Then $g_1=[d,f]=h$ and $h_1=[e,c]=g^{-1}$. Clearly, $g_1,h_1$ are in the center of $\R$. Also, it is easy to check that $g_1^{-1}=h^{-1}=d_1^{3^{n-1}}e_1^{3^{n-1}}f_1^{3^{n-1}}$,
and $h_1=g^{-1}=c_1^{3^{n-1}}d_1^{3^{n-1}}e_1^{3^{n-1}}$.

Finally, $[c_1,d_1]=[d,e]=d^{-3^{n-1}}e^{3^{n-1}}=c_1^{-3^{n-1}}d_1^{3^{n-1}}$,  $[c_1,f_1]=[d,c]=d^{-3^{n-1}}c^{3^{n-1}}=c_1^{-3^{n-1}}f_1^{3^{n-1}}$,
$[d_1,e_1]=[e,f]=e^{-3^{n-1}}f^{3^{n-1}}=d_1^{-3^{n-1}}e_1^{3^{n-1}}$, and
$[e_1,f_1]=[f,c]=f^{-3^{n-1}}c^{3^{n-1}}=e_1^{-3^{n-1}}f_1^{3^{n-1}}$.
This proves Claim~1.\medskip

\f{\bf Claim~2}\ $\R$ has no automorphisms swapping $a$ and $b$. \medskip

Suppose on the contrary that $\a$ is an automorphism of $\R$ such that $a^\a=b$ and $b^\a=a$.
Then $c^\a=[b,a]=c^{-1}$, $d^\a=[a,b^2]=f^{-1}$, $e^\a=[b^2,a^2]=e^{-1}$ and $f^\a=[a^2,b]=d^{-1}$.
Furthermore, $g^\a=[c^\a,e^\a]=[c^{-1},e^{-1}]=g^{e^{-1}c^{-1}}=g$ and $h^\a=[d^\a,f^\a]=[f^{-1},d^{-1}]=h^{-1}$.
As $g^{-1}=d^{3^{n-1}}e^{3^{n-1}}f^{3^{n-1}}$, one has $d^{3^{n-1}}e^{3^{n-1}}f^{3^{n-1}}=g^{-1}=(g^{-1})^\a=f^{-3^{n-1}}e^{-3^{n-1}}d^{-3^{n-1}}$.
This forces that $g^{-1}=d^{3^{n-1}}e^{3^{n-1}}f^{3^{n-1}}=1$, a contradiction. \medskip

Now we are ready to complete the proof. Note that $c^{3^{n-1}}, d^{3^{n-1}}, e^{3^{n-1}}$ are in the center of $D$, and $f^{3^{n-1}}=c^{3^{n-1}}d^{-3^{n-1}}e^{3^{n-1}}$. Set $M=\lg c^{3^{n-1}}, d^{3^{n-1}}, e^{3^{n-1}}\rg$. Then $M\cong C_3^3$.
By the following relations,
\[
\begin{array}{l}
h=c^{3^{n-1}}d^{3^{n-1}}e^{3^{n-1}}, g^{-1}=d^{3^{n-1}}e^{3^{n-1}}f^{3^{n-1}}, \\

[c,d]=c^{-3^{n-1}}d^{3^{n-1}}, [c,f]=c^{-3^{n-1}}f^{3^{n-1}},[d,e]=d^{-3^{n-1}}e^{3^{n-1}},[e,f]=e^{-3^{n-1}}f^{3^{n-1}}.
\end{array}
\]
we conclude that $M\leq D'$. It follows that $D/M=\lg cM\rg\times \lg dM\rg\times \lg eM\rg\times \lg fM\rg\cong C_{3^{n-1}}^4$.
Since $\R/D\cong C_3^2$ and $M\cong C_3^3$, it follows that $|\R|=3^{4n+1}$.

By Claim~1, the map $a\mapsto b, b\mapsto a^2$ induces an automorphism, say $\b$, of $\R$, and $\b$ cyclically permutates the elements in $S$.
So, $\b\in \Aut(\R,S)$. Since $S=\{a,a^2,b,b^2\}$, one has $\Aut(\R,S)\leq D_8$. By Claim~2, $\R$ has no automorphisms swapping $a$ and $b$. Consequently, we have $\Aut(\R,S)\cong C_4$. Since $\Upsilon$ is a tetravalent graph of order $3^{4n+1}$, $\Aut(\Upsilon)$ is a $\{2,3\}$-group, and so it is a solvable group. By Lemma~\ref{lem:testfornormal2geod-t}, $\Upsilon$ is $3$-CSH but not $2$-geodesic-transitive, and $C(\Upsilon)$ is a trivalent symmetric graph of type $2^2$.\hfill\qed

\f{\bf Remark on Lemma~\ref{lem:cubic2-2type}.}\  (1)\ We also verify Lemma~\ref{lem:cubic2-2type} in case $n=2$ by using {\tt Magma~{\rm\cite{BCP}}}, and the reader may see Appendix~8 for the {\tt Magma} programs.

(2)\ In 2006, Feng and Kwak~\cite{Feng-Kwak} posed the following conjecture.

\medskip
\f{\bf Conjecture}\ {\em Every connected trivalent symmetric graph of order $2\cdot 3^m$ is a Cayley graph for each $m\geq 1$. }
\medskip

By Lemma~\ref{lem:cubic2-2type}, $C(\Upsilon)$ is a trivalent symmetric graph of type $2^2$ and of order $2\cdot 3^{4n}$ with $n\geq 1$. If $C(\Upsilon)$ is a Cayley graph, then it would have an automorphism of order $2$ which swaps the two vertices of an edge. However, this is impossible since every edge-stabilizer for $C(\Upsilon)$ is isomorphic to $C_4$. Consequently, $C(\Upsilon)$ is a non-Cayley graph. This implies that the above conjecture is not true.

\section{Appendix: Magma programs in Section~6}

\f{\bf Appendix~1} (Programs for the graph $\Delta$ in Construction~I): First, we input a group
$G\lg a,b,c,d,e,f,g,h\rg:={\rm Group}\lg a,b,c,d,e,f,g,h\ |\ R \rg$, where $R$ is a set of relations as given in Construction~I.\smallskip

Construction of the group group $\H$: {\tt H,q:=pQuotient(G,2,100)};\smallskip

The order of group $\H$: {\tt FactoredOrder(H);}\smallskip

The derived subgroup of $\H$: {\tt D:=DerivedSubgroup(H);}\smallskip

The derived subgroup of $\H'$: {\tt DD:=DerivedSubgroup(D);}\smallskip

$\H/\H'\cong C_2^4$, $\H'/\H''\cong C_2^4\times C_4^2$ and $\H''\cong C_2$:

{\tt GroupName(H/D); GroupName(D/DD); GroupName(DD);}\smallskip

Construction of subgroups $A$ and $B$:

{\tt a:=a@q; b:=b@q; e:=e@q; f:=f@q; A:=sub<H|a,b >; B:=sub<H|e,f>;} \smallskip

Test $A\cong B\cong C_2^2$ and $A\cap B=1$:

{\tt $\sharp$A;
IsElementaryAbelian(A);

$\sharp$B;
IsElementaryAbelian(B);

$\sharp$(A meet B); H eq sub<H|a,b,e,f>;}
\smallskip

The following maps are not automorphisms of $\H$:

{\tt hom<H->H|a->b,b->a,e->f,f->e>;

hom<H->H|a->b,b->a,e->e,f->f>;

hom<H->H|a->a,b->b,e->f,f->e>;
}

The following maps are automorphisms of $\H$:

{\tt alpha:=hom<H->H|a->b,b->a*b,e->e,f->f>;
alpha;
Kernel(alpha);
Image(alpha) eq H;

beta:=hom<H->H|a->a,b->b,e->f,f->e*f>;
beta;
Kernel(beta);
Image(beta) eq H;

gamma:=hom<H->H|a->e,b->f,e->a,f->b>;
gamma;
Kernel(gamma);
Image(gamma) eq H;
}

\smallskip

Construction of subgroup $K$ and testing $K\cap A=K\cap B=1$:

{\tt K:=sub<H|a*e,b*f,DerivedSubgroup(H)>;

$\sharp$(K meet A);
$\sharp$(K meet B);}

\medskip
\f{\bf Appendix~2} (Programs for the construction of a Cayley graph):

{\tt Cay:=function(G,S);

V:={g:g in G};

E:={{g,s*g}:g in G,s in S};

return Graph<V|E>;

end function;}

\medskip
\f{\bf Appendix~3} (Programs for the graph $\Theta$ in Example~\ref{exam:5-val-1}): First, we input a group

{\tt G<a,b,c>:=Group<a,b,c| $a^5,b^5,c^5$, c=(a,b), (a,c)=(b,c)=1>;}\smallskip

Construction of the group $\M$: {\tt M,q:=pQuotient(G,5,100)};\smallskip

The order of group $\M$: {\tt FactoredOrder(M)};\smallskip

The derived subgroup of $\M$: {\tt D:=DerivedSubgroup(M);}\smallskip

The derived subgroup of $\M'$: {\tt DD:=DerivedSubgroup(D);}\smallskip

$\M/\M'\cong C_5^2$, $\M'/\M''\cong C_5$ and $\M''=1$:

{\tt GroupName(M/D); GroupName(D/DD); GroupName(DD);}\smallskip

Construction of $\Theta=\Cay(\M, S)$:

{\tt a:=a@q; b:=b@q;

S:=$\{a,a^2,a^3,a^4,b,b^2,b^3,b^4\}$;

Theta:=Cay(M,S);}

\smallskip
Automorphism group of $\Theta$:

{\tt A:=AutomorphismGroup(Theta);}

\smallskip
$\Theta$ is a normal Cayley graph on $\M$ (We find that every Sylow $5$-subgroup of $\Aut(\Theta)$ is normal and regular on $V(\Theta)$. This implies that $\Theta$ is normal):

{\tt P:=SylowSubgroup(A,5);

IsNormal(A,P);

IsRegular(P);}

\smallskip
$\Aut(\M,S)\cong C_4\wr C_2$:

{\tt A1:=Stabilizer(A,1);

GroupName(A1);}

\medskip
\f{\bf Appendix~4} (Programs for the graph $\Phi$ in Example~\ref{exam:5-val-2}): First, we input a group

{\tt G<a,b,c,d,e>:=Group<a,b,c,d,e| $a^5,b^5,c^5,d^5,e^5$, c=(a,b),d=(a,c),e=(b,c),

(a,d)=(b,d)=(a,e)=(b,e)=1>;}\smallskip

Construction of the group group $\T$: {\tt T,q:=pQuotient(G,5,100)};\smallskip

The order of group $\T$: {\tt FactoredOrder(T)};\smallskip

The derived subgroup of $\T$: {\tt D:=DerivedSubgroup(T);}\smallskip

The derived subgroup of $\T'$: {\tt DD:=DerivedSubgroup(T);}\smallskip

$\T/\T'\cong C_5^2$, $\T'/\T''\cong C_5^3$ and $\T''=1$:

{\tt GroupName(T/D); GroupName(D/DD); GroupName(DD);}\smallskip

Construction of $\Phi=\Cay(\T, S)$:

{\tt a:=a@q; b:=b@q;

S:=$\{a,a^2,a^3,a^4,b,b^2,b^3,b^4\}$;

Phi:=Cay(T,S);}

\smallskip
Automorphism group of $\Phi$:

{\tt A:=AutomorphismGroup(Phi);}

\smallskip
$\Phi$ is a normal Cayley graph on $\T$ (We find that every Sylow $5$-subgroup of $\Aut(\Phi)$ is normal and regular on $V(\Phi)$. This implies that $\Phi$ is normal):

{\tt P:=SylowSubgroup(A,5);

IsNormal(A,P);

IsRegular(P);}

\smallskip
$\Aut(\T,S)\cong C_4\wr C_2$:

{\tt A1:=Stabilizer(A,1);

GroupName(A1);}

\medskip
\f{\bf Appendix~5} (Programs for the graph $\Theta$ in Example~\ref{exam:cond-3}): First, we input a group

{\tt G<a,b,c,d,e,f,g,x,y,z>:=Group<a,b,c,d,e,f,g,x,y,z | $a^2,b^2,c^2,e^2,f^2,g^2, x^2,y^2,z^2$,

(a,b)=(a,c)=(b,c)=1, (e,f)=(e,g)=(f,g)=1, (e,a)=x*y*z, (e,b)=x*z, (e,c)=x,

(f,a)=x*z, (f,b)=x, (f,c)=y, (g,a)=x, (g,b)=y, (g,c)=z >;}
\smallskip

Construction of the group group $\mathcal{G}$: {\tt G,q:=pQuotient(G,2,100)};\smallskip

The order of group $\mathcal{G}$: {\tt FactoredOrder(G)};\smallskip

The derived subgroup of $\mathcal{G}$: {\tt D:=DerivedSubgroup(G);}\smallskip

The derived subgroup of $\mathcal{G}'$: {\tt DD:=DerivedSubgroup(G);}\smallskip

$\mathcal{G}/\mathcal{G}'\cong C_2^6$, $\mathcal{G}'/\mathcal{G}''\cong C_2^3$ and $\mathcal{G}''=1$:

{\tt GroupName(G/D); GroupName(D/DD); GroupName(DD);}\smallskip

Construction of subgroups $A$ and $B$, and testing $A\cong B\cong C_2^3$, $A\cap B=1$ and $\mathcal{G}=\lg A,B\rg$:

{\tt a:=a@q;
b:=b@q;
c:=c@q;
e:=e@q;
f:=f@q;
g:=g@q;

A:=sub<G|a,b,c>;
$\sharp$A;
IsElementaryAbelian(A);

B:=sub<G|e,f,g>;
$\sharp$B;
IsElementaryAbelian(B);

$\sharp$(A meet B);

H eq sub<G|A,B>;}

\smallskip
Construction of $\Theta=\Cay(\mathcal{G}, S)$:

{\tt
S:=$\{$x:x in A|x ne G!1$\}$ join $\{$y:y in B|y ne G!1$\}$;

Theta:=Cay(G,S);}

\smallskip
Automorphism group of $\Theta$:

{\tt au:=AutomorphismGroup(Theta);}

\smallskip
$\Theta$ is a normal Cayley graph on $\mathcal{G}$ (Note that $R(\mathcal{G})$ is a regular subgroup of $\Aut(\Theta)$. We first list all regular subgroups of $\Aut(\Theta)$, and then we find that among these subgroups, there is only one which is isomorphic to $R(\mathcal{G})$ and this subgroup is normal in $\Aut(\Theta)$):

{\tt R:=RegularSubgroups(au);//Find all regular subgroups of $\Aut(\Theta)$

T:=$\{\}$;

for i in $\{$1..$\sharp$R$\}$ do

if $\sharp$(DerivedSubgroup(G)) eq  $\sharp$(DerivedSubgroup(R[i]`subgroup)) then

Include($^\sim$T,i);

end if;

end for;

$\sharp$T;//$|T|=1$

i:=Random(T);

IsNormal(au,R[i]`subgroup);}

\smallskip
$\Aut(\mathcal{G},S)\cong C_7\times C_7)\rtimes C_6$:

{\tt au1:=Stabilizer(au,1);

GroupName(au1);}

\medskip
\f{\bf Appendix~6} (Programs for the graph $\Pi$ in Example~\ref{exam:cond-4}): First, we input a group
$G\lg a,b,c,d,e,u,v,x,y,z,f,g,h,i,j\rg:={\rm Group}\lg a,b,c,d,e,u,v,x,y,z,f,g,h,i,j\ |\ R \rg$, where $R$ is a set of relations as given in Example~\ref{exam:cond-4}.\smallskip

Construction of the group group $\L$: {\tt L,q:=pQuotient(G,2,100)};\smallskip

The order of group $\L$: {\tt FactoredOrder(L)};\smallskip

The derived subgroup of $\L$: {\tt D:=DerivedSubgroup(L);}\smallskip

The derived subgroup of $\L'$: {\tt DD:=DerivedSubgroup(L);}\smallskip

$\L/\L'\cong C_2^{10}$, $\L'/\L''\cong C_2^5$ and $\L''=1$:

{\tt GroupName(L/D); GroupName(D/DD); GroupName(DD);}

Construction of subgroups $A$ and $B$, and testing $A\cong B\cong C_2^5$, $A\cap B=1$ and $\L=\lg A,B\rg$:

{\tt a:=a@q; b:=b@q; c:=c@q; d:=d@q; e:=e@q;

u:=u@q;
v:=v@q;
x:=x@q;
y:=y@q;
z:=z@q;

A:=sub<L|a,b,c,d,e>;
$\sharp$A;
IsElementaryAbelian(A);

B:=sub<L|u,v,x,y,z>;
$\sharp$B;
IsElementaryAbelian(B);

$\sharp$(A meet B);
H eq sub<L|A,B>;
}

\smallskip
Construction of the clique graph of $\Pi=\Cay(\L, S)$:

{\tt
V1:=$\{\}$;

for h in L do

Vh:=$\{\}$;

for w in A do

Include($^\sim$Vh,w*h);

end for;

Include($^\sim$V1,Vh);

end for;

V2:=$\{\}$;

for h in L do

Vh:=$\{\}$;

for w in B do

Include($^\sim$Vh,w*h);

end for;

Include($^\sim$V2,Vh);

end for;

V:=V1 join V2;

E:=$\{\}$;

for w1 in V1 do

for w2 in V2 do

if $\sharp$(w1 meet w2) eq 1 then

Include($^\sim$E,$\{$w1,w2$\}$);

end if;

end for;

end for;

CPi:=Graph<V|E>;//This is the clique graph of $\Pi$
}

\smallskip
Automorphism group of $C(\Pi)$:

{\tt au:=AutomorphismGroup(CPi);}

\smallskip
$\Pi$ is a normal Cayley graph on $\L$ (Note that $R(\L)$ is a subgroup of $\Aut(C(\Pi))$. We first list all subgroups of $\Aut(C(\Pi))$ of order $2^{15}$, and we find that among these subgroups, there is only one which is isomorphic to $R(\L)$ and this subgroup is normal in $\Aut(C(\Pi))$):

{\tt R:=Subgroups(au:OrderEqual:=$2^{15}$);

T:=$\{\}$;

for i in $\{$1..$\sharp$R$\}$ do

if $\sharp$(DerivedSubgroup(L)) eq  $\sharp$(DerivedSubgroup(R[i]`subgroup)) then

Include($^\sim$T,i);

end if;

end for;

$\sharp$T;//$|T|=1$

i:=Random(T);

IsNormal(au,R[i]`subgroup);}

\smallskip
$\Aut(\R,S)\cong (C_{31}\times C_{31})\rtimes C_{10}$ (Note that $\Aut(\R,S)\cong \Aut(C(\Pi))/R(\L)$):

{\tt GroupName(au/R[i]`subgroup);}

\medskip
\f{\bf Appendix~7} (Programs for the graph $\Lambda$ in Example~\ref{exam:val-5-m16}): First, we input a group

{\tt G<a,b,c,d,e,f,g,h,k>:=Group<a,b,c,d,e,f,g,h,k| $a^5,b^5,c^5,d^5,e^5,f^5,g^5,h^5,k^5,$

c=(a,b),d=(a,c),e=(b,c),(d,e)=1, (a,d)=f, (b,d)=g, (a,e)=h, (b,e)=k,

(a,f)=(a,g)=(a,h)=(a,k)=(b,f)=(b,g)=(b,h)=(b,k)=1, f=$k^{-2}$, g=$h^{-2}$ >;
}\smallskip

Construction of the group group $\N$: {\tt N,q:=pQuotient(G,5,100)};\smallskip

The order of group $\N$: {\tt FactoredOrder(N)};\smallskip

The derived subgroup of $\N$: {\tt D:=DerivedSubgroup(N);}\smallskip

The derived subgroup of $\N'$: {\tt DD:=DerivedSubgroup(N);}\smallskip

$\N/\N'\cong C_5^2$, $\N'/\N''\cong C_5^4$ and $\N''=1$:

{\tt GroupName(N/D); GroupName(D/DD); GroupName(DD);}\smallskip

Construction of $\Lambda=\Cay(\N, S)$:

{\tt a:=a@q; b:=b@q;

S:=$\{a,a^2,a^3,a^4,b,b^2,b^3,b^4\}$;

Lam:=Cay(N,S);}

\smallskip
Automorphism group of $\Lambda$:

{\tt A:=AutomorphismGroup(Lam);}

\smallskip
$\Lambda$ is a normal Cayley graph on $\N$ (We find that every Sylow $5$-subgroup of $\Aut(\Lambda)$ is normal and regular on $V(\Lambda)$. This implies that $\Lambda$ is normal):

{\tt P:=SylowSubgroup(A,5);

IsNormal(A,P);

IsRegular(P);}

\smallskip
$\Aut(\N,S)\cong {\rm M}_{16}$ (Note that ${\rm M}_{16}$ is just the semidihedral group of order $16$):

{\tt A1:=Stabilizer(A,1);

GroupName(A1);}

\medskip
\f{\bf Appendix~8} (Programs for the graph $\Upsilon$ in Construction~II (in case $n=2$)): First, we input a group

{\tt G<a,b,c,d,e,f,g,h>:=Group<a,b,c,d,e,f,g,h | $a^3,b^3,c^ 9,d^9,e^9,f^9,g^3,h^3,$

c=(a,b),d=$(b,a^2)$,e=$(a^2,b^2)$,f=$(b^2,a)$, (c,d)=$c^{-3}*d^3$,(c,f)=$c^{-3}*f^3$,(d,e)=$d^{-3}*e^3$,

(e,f)=$e^{-3}*f^3$, $f^3=c^3*d^{-3}*e^3$,g=(c,e),h=(d,f), h=$c^3*d^3*e^3$, $g^{-1}=d^3*e^3*f^3$,

$(c^3,d)=(c^3,e)=(c^3,f)=1,(d^3,c)=(d^3,e)=(d^3,f)=1,$

$(e^3,c)=(e^3,d)=(e^3,f)=1$>;
}\smallskip

Construction of the group group $\R$: {\tt R,q:=pQuotient(G,3,100)};\smallskip

The order of group $\R$: {\tt FactoredOrder(R)};\smallskip

The derived subgroup of $\R$: {\tt D:=DerivedSubgroup(R);}\smallskip

The derived subgroup of $\R'$: {\tt DD:=DerivedSubgroup(R);}\smallskip

$\R/\R'\cong C_3^5$, $\R'/\R''\cong C_3^3\times C_9$ and $\R''\cong C_3^2$:

{\tt GroupName(R/D); GroupName(D/DD); GroupName(DD);}\smallskip

Construction of $\Upsilon=\Cay(\R, S)$:

{\tt a:=a@q; b:=b@q;

S:=$\{a,a^2,b,b^2\}$;

Ups:=Cay(R,S);}

\smallskip
Automorphism group of $\Upsilon$:

{\tt A:=AutomorphismGroup(Ups);}

\smallskip
$\Upsilon$ is a normal Cayley graph on $\R$ (We find that every Sylow $3$-subgroup of $\Aut(\Upsilon)$ is normal and regular on $V(\Upsilon)$. This implies that $\Upsilon$ is normal):

{\tt P:=SylowSubgroup(A,3);

IsNormal(A,P);

IsRegular(P);}

\smallskip
$\Aut(\R,S)\cong C_{4}$:

{\tt A1:=Stabilizer(A,1);

GroupName(A1);}

\medskip
\f {\bf Acknowledgements:}\
This work was supported by the National Natural Science Foundation of China (12071023).

{}
\end{document}